\documentclass[preprint,12pt]{elsarticle}

\usepackage{algorithm}
\usepackage{algorithmic}
\usepackage{arydshln}

\usepackage{amsmath}
\usepackage{amssymb}

\begin{document}

\begin{frontmatter}




\title{Domain Decomposition Subspace Neural Network Method for Solving Linear and Nonlinear Partial Differential Equations}%

\author[1]{Zhenxing Fu}
\cortext[cor1]{Corresponding author}
\author[1]{Hongliang Liu} 
\author[2,3]{Zhiqiang Sheng\corref{cor1}}
\ead{sheng_zhiqiang@iapcm.ac.cn}
\author[1]{Baixue Xing}

\address[1]{School of Mathematics and Computational Science, Xiangtan University, Hunan Key Laboratory for Computational Science and Numerical Simulation, Xiangtan, Hunan, 411105, China.}
\address[2]{National Key Laboratory of Computational Physics, Institute of Applied Physics and Computational Mathematics, Beijing, 100088, China.}
\address[3]{HEDPS, Center for Applied Physics and Technology, and College of Engineering, Peking University, Beijing, 100871, China.}

\begin{abstract}
This paper proposes a domain decomposition subspace neural network method for efficiently solving linear and nonlinear partial differential equations. 
By combining the principles of domain decomposition and subspace neural networks, the method constructs basis functions using neural networks to approximate PDE solutions. 
It imposes $C^k$ continuity conditions at the interface of subdomains, ensuring smoothness across the global solution.
Nonlinear PDEs are solved using Picard and Newton iterations, analogous to classical methods. 
Numerical experiments demonstrate that our method achieves exceptionally high accuracy, with errors reaching up to $10^{-13}$, while significantly reducing computational costs compared to existing approaches, including PINNs, DGM, DRM. The results highlight the method's superior accuracy and training efficiency.
\end{abstract}


\end{frontmatter}

\date{}



\section{Introduction}
Partial differential equations (PDEs) \cite{Evans2022} play a central role in science and engineering, widely used to describe complex phenomena such as fluid dynamics, heat conduction, and quantum mechanics. 
Traditional numerical methods like finite difference \cite{FDM2020}, finite element \cite{FEM2013}, and spectral methods \cite{Spectral2001} have established solid foundations for PDE solutions, but face challenges in high-dimensional and nonlinear scenarios due to mesh dependency and computational complexity \cite{Chen2022RandomFeature}.

In recent years, with the rapid development of machine learning, neural networks have emerged as promising tools for solving PDEs due to their powerful function approximation capabilities. Among them, Physics-Informed Neural Networks (PINNs) \cite{Raissi2019PINNs} emerged as a paradigm shift by embedding physical constraints into loss functions, enabling mesh-free solutions without labeled data \cite{Yang2022MO-PINNs}. 
However, PINNs struggle with convergence instability and gradient vanishing in high-frequency or complex boundary problems \cite{Basir2022CriticalPINNs,Nguwi2024DeepBranchingSolver,Zuo2023ResidualNetworks}, particularly when handling discontinuous interfaces or multiscale dynamics \cite{Wang2024NASPINNs,Siegel2023GreedyTraining}. 
Recent advances in adaptive architectures \cite{Lin2024AdaptiveNN} and tensor networks \cite{Wang2024MultiEigenpairs} partially mitigate these issues.

To enhance PINNs' adaptability, hybrid strategies have been proposed. Physics-Informed Radial Basis Networks \cite{Bai2023PIRBN} leverage local approximation for high-frequency problems, while boundary-fitted PINNs \cite{Xie2023BoundaryFittingPINNs} improve geometric flexibility at the cost of manual trial function design. 
Incremental learning frameworks like iPINNs \cite{Dekhovich2025iPINNs} and non-gradient optimizers \cite{Peng2023NonGradientMethod} address training instability. 
Parallel efforts in operator learning \cite{Li2024LocalNeuralOperator} and foundation models \cite{Ye2024PDEformer} demonstrate potential.

Domain decomposition methods (DDMs) offer a promising pathway by combining localized neural solvers with global continuity conditions. 
Early attempts in overlapping Schwarz methods \cite{Jagtap2020XPINNs} faced accuracy loss from interface error propagation, whereas non-overlapping DDMs struggled with flux estimation \cite{Liu2023cvPINNs}. 
Recent breakthroughs in compensated deep Ritz methods \cite{Shang2023RNNPG} and variational formulations enable reliable flux transmission across subdomains, enhancing parallelism and representation capacity. 
Similar decomposition strategies have shown promise in electromagnetic simulations \cite{Cobb2023MaxwellsEquation}.

Additionally, Local Extreme Learning Machines (LocELM) \cite{Dong2021LocELM} building upon classical ELM frameworks \cite{Huang2006ELM} and subspace method based on neural networks (SNN) \cite{Sheng2024SNN} have emerged as efficient methods for solving PDEs. 
Recent advances in randomized neural networks with Petrov-Galerkin methods \cite{Shang2023RNNPG} further demonstrate the potential of subspace approaches. 
LocELM builds localized models to solve both linear and nonlinear equations efficiently, but its sensitivity to randomly initialized weights may lead to convergence instability in highly nonlinear problems. 
On the other hand, SNN approach the problem by constructing neural networks in subspaces, offering high solution accuracy at low computational cost.

Building upon these advances, this paper proposes a Domain Decomposition Subspace Neural Network (DD-SNN) method. 
By independently training local neural networks in different subdomains and imposing $C^k$continuity conditions at the subdomain interfaces, the proposed method ensures the smoothness of the solution, thus improving the solution accuracy while maintaining computational efficiency. 
Additionally, by optimizing the coefficients of basis functions through improved iteration schemes compared to conventional nonlinear least squares \cite{Peng2023NonGradientMethod}, the accuracy is further enhanced. Numerical experiments show that DD-SNN outperforms existing methods such as PINNs, LocELM, and SNN in both accuracy and computational efficiency.

The main contributions of this paper are as follows:

1. Combination of Domain Decomposition and Subspace Neural Networks: 
By dividing the computational domain into multiple subdomains and independently constructing local neural networks, the computational complexity of the global problem is significantly reduced. 
The subdomains are seamlessly connected through $C^k$ continuity conditions, resulting in a globally smooth solution, and parallel computation is supported, greatly enhancing the efficiency of solving complex problems.
2. Separated Training Strategy: 
This approach separates the training of basis functions from the solution approximation process, allowing different methods to optimize the basis functions, thus improving the flexibility and accuracy of the solution.
3. Efficient Solution of Nonlinear Problems: 
By using Picard and Newton iterations to optimize the coefficients of the basis functions, the method avoids the high computational cost of nonlinear least squares, resulting in a significant improvement in accuracy while reducing the solution time.
4. High Accuracy with Low Cost: 
In several numerical experiments, the error reaches the level of $10^{-13}$, outperforming existing methods such as ELM, DGM, and PINNs.

Extensive numerical experiments are conducted to evaluate the performance of the DD-SNN method in solving linear and nonlinear, stationary and time-dependent PDEs. 
The results show that DD-SNN significantly outperforms other methods in both accuracy and computational efficiency, with much lower numerical errors and training times.

The remainder of this paper is organized as follows: Section 2 introduces the DD-SNN method, including its principles, domain decomposition strategy, and local subspace neural network construction. 
Sections 3 and 4 detail the solution process for linear and nonlinear PDEs, respectively. Section 5 presents numerical experiments, comparing DD-SNN with PINNs, LocELM, SNN, and other methods in terms of accuracy and computational efficiency. 
Finally, Section 6 concludes the paper and discusses the potential and future developments of the DD-SNN method in solving complex PDEs.

\section{Domain Decomposition Subspace Neural Network Method}
Consider the partial differential equation
\begin{subequations}
    \begin{alignat}{2}
        \mathcal{L}u(\boldsymbol{x}, t) +  \mathcal{N}(u, u_t, \nabla u,\nabla^2 u, \cdots) & = f(\boldsymbol{x}, t), &\quad & (\boldsymbol{x}, t) \in \Omega,  \label{eq:1a} \\
        u(\boldsymbol{x}, t) & = g(\boldsymbol{x}, t), &\quad & (\boldsymbol{x}, t) \in \partial \Omega,  \label{eq:1b} 
    \end{alignat}
    \label{eq:1}
\end{subequations}
where $\boldsymbol{x} = (x_1, x_2, \ldots, x_d)^T$, $\Omega$ is a bounded domain in $\mathbb{R}^{d+1}$, $(\boldsymbol{x},t) \in \Omega \subset \mathbb{R}^d \times \mathbb{R}$, and when time is not involved, $\Omega$ represents a bounded domain in $\mathbb{R}^d$. The boundary $\partial \Omega$ represents the spatial or spatiotemporal boundary. 
The operator $\mathcal{L}$ is a linear differential operator with respect to $(\boldsymbol{x}, t)$, $\mathcal{N}$ is a nonlinear operator, $f(\boldsymbol{x}, t)$ is a known source term, and $g(\boldsymbol{x}, t)$ represents boundary conditions.

The domain $\Omega$ is divided into $N_k$ non-overlapping subdomains as follows:
$$\Omega=\bigcup_{k=1}^{N_k}\Omega_k, \quad  \Omega_p \cap \Omega_q = \emptyset,\quad \forall p \neq q,$$
where each subdomain $\Omega_k$ corresponds to a local solving region.

\begin{figure}[htbp]
  \centering
  \includegraphics[width=0.7\textwidth]{ 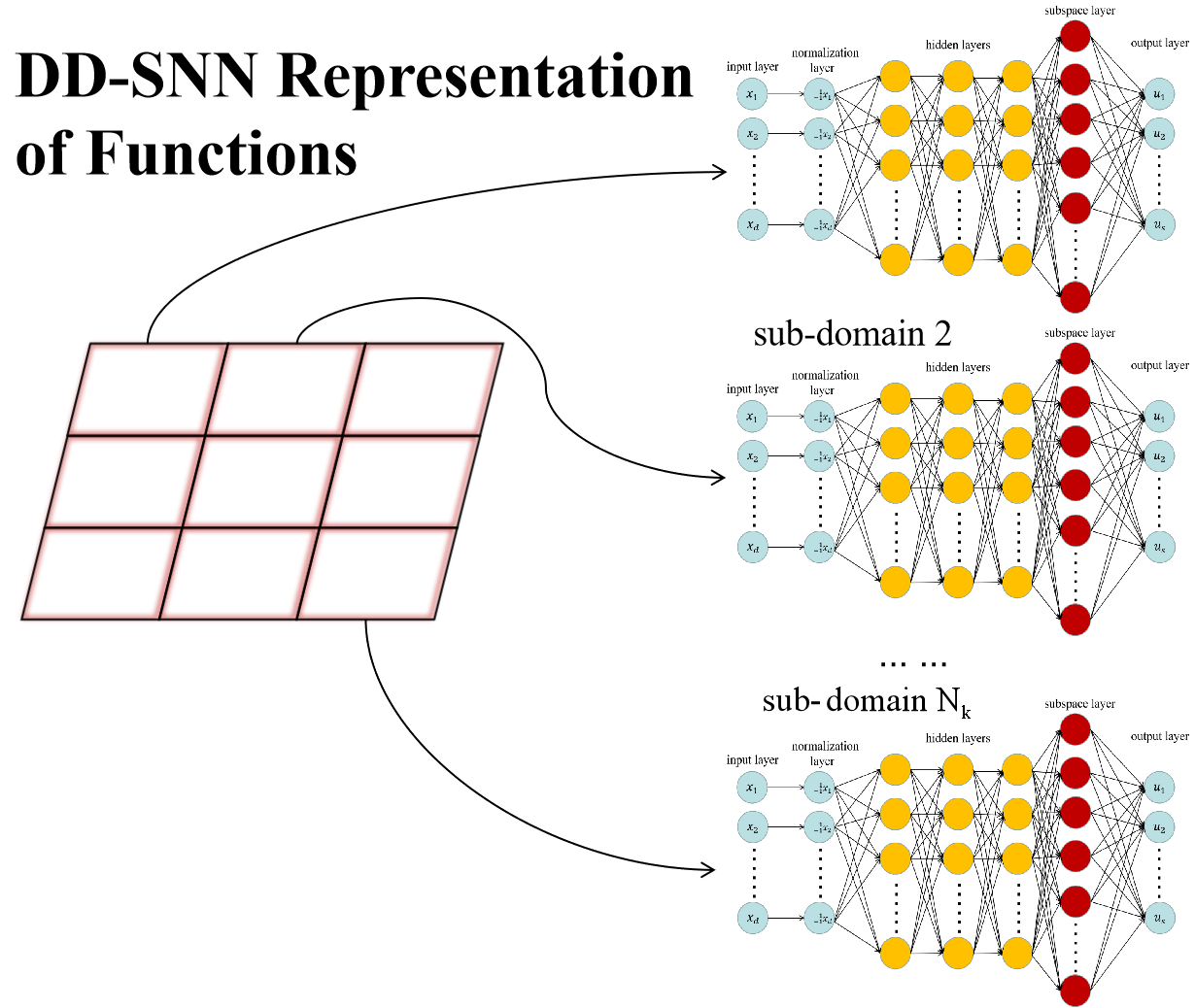}
  \caption{Neural network framework.}
  \label{fig:1}
\end{figure}

Within each subdomain $\Omega_k$, the local solution $u^k(\boldsymbol{x},t)$ is approximated using a subspace neural network that consists of an input layer, a normalization layer, hidden layers, a subspace layer, and an output layer~\cite{Sheng2024SNN}, as shown in Fig.~\ref{fig:1}. 
Let the total number of hidden layers be $L$, with the number of neurons in each hidden layer given by $n^{(1)}, n^{(2)}, \dots, n^{(L)}$. 
The subspace layer has dimension $M$, and its basis functions are denoted as $\Phi_j (j = 1, 2, \dots, M)$, which span a subspace that approximates the PDE solution. The output layer has dimension $s$, corresponding to the number of solution components.

The propagation process of the SNN in subdomain $\Omega_k$ is given as
\begin{subequations}
  \begin{align}
      &\mathbf{z}_0^k = X_{norm}^k(\boldsymbol{x},t) = \begin{bmatrix} \frac{2}{b_{\boldsymbol{x}}^k - a_{\boldsymbol{x}}^k} \cdot (\boldsymbol{x} - a_{\boldsymbol{x}}^k) - 1 \\
      \frac{2}{b_t^k - a_t^k} \cdot (t - a_t^k) - 1 \end{bmatrix},  \label{eq:2a} \\
      &\mathbf{z}_l^k = \sigma^k_l(W_l^k \cdot \mathbf{z}_{l-1}^k + b_l^k), \quad l = 1, 2, \dots, L, \label{eq:2b} \\
      &\mathbf{z}_{L+1}^k = \Phi^k(\boldsymbol{x}, t) = (\phi_1^k, \phi_2^k, \dots, \phi_M^k)^{T},  \label{eq:2c} \\
      &\hat{u}^k(\boldsymbol{x}, t) = \left[\hat{u}^k_1, \hat{u}^k_2,\dots,\hat{u}^k_s\right] = \sum^{M}_{j=1} \beta_{ji}^k \cdot \phi_j^k, \quad 1 \leq i \leq s, \label{eq:2d}
  \end{align}
\end{subequations}
where $X_{norm}^k(\boldsymbol{x},t)$ is a normalization function mapping the input $(\boldsymbol{x},t)$ to the interval $[-1,1]^{d+1}$. 
The parameters $a^k_{\boldsymbol{x}}, b^k_{\boldsymbol{x}}$ are the lower and upper bounds of the spatial direction in subdomain $\Omega_k$, while $a^k_t, b^k_t$ are the lower and upper bounds of the temporal direction. 
$W_l^k \in \mathbb{R}^{n^{(l)} \times n^{(l-1)}}$ and $b_l^k \in \mathbb{R}^{n^{(l)}}$ are the weight matrix and bias vector, respectively, for subdomain $\Omega_k$. 
The activation function is denoted as $\sigma_l^k(\cdot)$, and $\Phi^k(\boldsymbol{x}, t)$ represents the basis functions in the subdomain. 
The parameter set $\theta^k = \{W_1^k,W_2^k, \dots, W_{L}^k; b_1^k,b_2^k, \dots, b_{L}^k\}$, along with the weight vector $\boldsymbol{\beta}^k$, determines the network output $u^k(\boldsymbol{x},t; \theta^k, \boldsymbol{\beta}^k)$, denoted as $\hat{u}^k(\boldsymbol{x},t)$.

The network parameters $\theta^k$ and the subspace basis functions are optimized by minimizing a local loss function $\mathcal{L}_{PDE}^{k}(\boldsymbol{x},t;\theta^k,\boldsymbol{\beta}^k)$ that only contains information about the PDE: 
\begin{eqnarray}
    \mathcal{L}_{PDE}^{k}(\boldsymbol{x}, t;\theta^k, \boldsymbol{\beta}^k)
     &= &\frac{1}{N_r} \sum_{i=1}^{N_r} \left( \mathcal{L}\hat{u}^k(\boldsymbol{x}^i, t^i) +  \mathcal{N}(\hat{u}^k(\boldsymbol{x}^i,t^i),\hat{u}^k_t(\boldsymbol{x}^i, t^i),\right.\nonumber\\ 
     &&\left.\nabla \hat{u}^k(\boldsymbol{x}^i, t^i), \nabla^2 \hat{u}^k(\boldsymbol{x}^i,t^i), \dots) - f^k(\boldsymbol{x}^i, t^i) \right)^2,
    \label{eq:3}
  \end{eqnarray}
where $N_r$ is the number of configuration points. 
During training, the network parameters $\theta^k$ are updated by minimizing the loss function \eqref{eq:3}.  
Although the basis combination coefficients $\boldsymbol{\beta}^k$ can also be updated during training, their final values are determined by enforcing the equations and boundary conditions at collocation points. Therefore, $\boldsymbol{\beta}^k$ is set to 1 during training to ensure balanced contributions from each basis function.

To balance computational efficiency and solution accuracy, training stops when either the relative reduction in the loss function \eqref{eq:3} with respect to its initial value $\mathcal{L}_{PDE_0}^k(\boldsymbol{x}, t;\theta^k,\boldsymbol{\beta}^k)$ satisfies
\begin{equation}
  \frac{\mathcal{L}_{PDE}^k(\boldsymbol{x}, t;\theta^k,\boldsymbol{\beta}^k)}{\mathcal{L}_{PDE_0}^k(\boldsymbol{x}, t;\theta^k,\boldsymbol{\beta}^k)} \leq \varepsilon,
\end{equation}
or the maximum number of training epochs $N_{max}$ is reached. 
Here, $\varepsilon$ is a small positive constant, typically set to $10^{-3}$.

Remark 2.1: The training processes on different subdomains are completely independent, inherently allowing for parallel computation, thereby significantly improving computational efficiency.

If two subdomains $\Omega_p$ and $\Omega_q$ share a boundary, the shared boundary is denoted as $\mathcal{I}_{pq}$. 
On these shared interfaces, we impose continuity conditions:
\begin{equation}
  \frac{\partial^\alpha \hat{u}^p(\boldsymbol{x},t)}{\partial x^\alpha_s}-\frac{\partial^\alpha \hat{u}^q(\boldsymbol{x},t)}{\partial x^\alpha_s}=0, \quad \alpha=1,2,\dots,k_s, \quad \forall (\boldsymbol{x},t) \in \mathcal{I}_{pq}, \label{eq:5}
\end{equation}
where $x_s$ represents the $s$-th coordinate direction, $\alpha$ denotes the order of derivatives $(0 \leq \alpha \leq k_s)$, and $k_s$ is the continuity order along the $x_s$ direction. 
This ensures continuity of the function values and their partial derivatives up to order $k_s$, thereby maintaining the smoothness of the global solution. 
The choice of continuity order $k_s$ is closely related to the properties of the equation. 
For example, if the highest-order derivative in a specific direction is $k_s$, the solution is usually required to have $C^{k_s-1}$ continuity across the subdomain interfaces in that direction.

After training the basis functions, the local solution $\hat{u}^k(\boldsymbol{x},t)$ in each subdomain $\Omega_k (1\leq k \leq N_k)$ is provided by Eq.~\eqref{eq:2d}. 
Subsequently, continuity conditions are imposed on the shared interfaces $\mathcal{I}_{pq}$ to construct a global algebraic system. Depending on whether the equation is linear or nonlinear, appropriate methods such as least squares, Picard iteration, or Newton iteration are employed to solve for the optimal combination coefficients $\boldsymbol{\beta}^k$ and update these coefficients. 
Finally, the local solutions from all subdomains are assembled in the global subspace to obtain the approximate solution $\hat{u}(\boldsymbol{x},t)$. The algorithm flow is shown in Algorithm~\ref{alg1:DD-SNN}.

It is worth noting that when the number of subdomains $N_k$ is reduced to 1, such that the entire domain $\Omega$ is covered by a single neural network, the method simplifies to a standard SNN. 
Furthermore, if both the number of hidden layers and the number of epochs are reduced to zero, the method degenerates into an extreme learning machine (ELM). 
Additionally, if the number of epochs are reduced to zero while keeping the randomly initialized weights and biases fixed within the range $[-R_m, R_m]$, the method reduces to a local extreme learning machine (LocELM).

\begin{algorithm}[htbp]
  \caption{DD-SNN Algorithm Flow}
  \label{alg1:DD-SNN}
  \begin{algorithmic}[1]
      \STATE \textbf{Step 1}: Partition the domain $\Omega$ into $N_k$ independent, non-overlapping subdomains $\{\Omega_k\}_{k=1}^{N_k}$.
      \STATE \textbf{Step 2}: Construct the loss function $\mathcal{L}_{PDE}^k$ for each subdomain and independently minimize these loss functions to update neural network parameters.
      \STATE \textbf{Step 3}: Obtain the basis functions in the subspace layer and represent the solution as a linear combination of these basis functions.
      \STATE \textbf{Step 4}: Enforce the equation within each subdomain, the boundary conditions on the domain interfaces, and the continuity conditions on the shared interfaces between subdomains. Solve the resulting algebraic system to determine the combination coefficients of the basis functions.
      \STATE \textbf{Step 5}: Assemble the local solutions to obtain the approximate global solution.
  \end{algorithmic}
\end{algorithm}

In the DD-SNN method, the optimization of network parameters $\theta^k$ and the solution of the weight vector $\boldsymbol{\beta}^k$ are two independent processes. 
This design allows the basis functions to be trained using different strategies and sampling methods. 
Additionally, the DD-SNN framework is highly flexible, allowing the local neural networks in different subdomains to adopt different hyperparameter settings (such as the depth and width of hidden layers, types of activation functions, and the dimension of the subspace) to better capture the characteristics and fine structures of local solutions in complex regions. 
However, to simplify implementation and ensure overall consistency and coherence, this study adopts identical hyperparameter configurations for all local neural networks across the subdomains. 
Our discussion below is divided into two cases: linear problem and nonlinear problem.

\section{Linear Equations}
For solving linear partial differential equations, the nonlinear term in Eq.~\eqref{eq:1a} is set to zero, i.e. ,
$\mathcal{N}(u, u_t,\nabla u,\nabla^2 u, \dots)=0$, resulting in the linear equation over the spatiotemporal domain $\Omega$:
\begin{subequations}
    \begin{alignat}{2}
        \mathcal{L}u(\boldsymbol{x}, t) & = f(\boldsymbol{x}, t), &\quad &(\boldsymbol{x},t) \in \Omega,  \label{eq:6a} \\
        u(\boldsymbol{x}, t) & = g(\boldsymbol{x}, t), &\quad &(\boldsymbol{x}, t) \in \partial \Omega.  \label{eq:6b}
    \end{alignat}
    \label{eq:6}
\end{subequations}

To solve Eq.~\eqref{eq:6}, following Steps 1 and 2 of Algorithm~\ref{alg1:DD-SNN}, the loss function is minimized within each subdomain $\Omega_k$ to optimize the parameters $\theta^k$, yielding a set of basis functions $\phi_j^k(j=1,2,\dots,M)$. 
The approximate solution is then expressed as a linear combination of these basis functions:
\begin{equation}
    \hat{u}^k(\boldsymbol{x},t) = \sum^{M}_{j=1}\beta_{ji}^k \cdot \phi_j^k(\boldsymbol{x},t).
    \label{eq:777}
\end{equation}

Substituting this solution into Eq.~\eqref{eq:6a} yields:
\begin{equation}
    \mathcal{L}\left(\sum^{M}_{j=1}\boldsymbol{\beta}^k_{ji}\cdot\phi_j^k(\boldsymbol{x},t) \right)=f^k(\boldsymbol{x},t) \quad 1\leq i \leq s, \quad (\boldsymbol{x},t) \in \Omega_k.
    \label{eq:7}
\end{equation}
Enforcing this equation at specific collocation points results in the following matrix form:
\begin{equation}
    \boldsymbol{A}^k \boldsymbol{\beta}^k = f^k,
    \label{eq:8}
\end{equation}
where $\boldsymbol{A}^k$ is the coefficient matrix representing the linear operator $\mathcal{L}$ acting on the basis functions in subdomain $\Omega_k$, $\boldsymbol{\beta}^k$ is the vector of combination coefficients for the basis functions, and $f^k$ is the discretized vector of the source term $f^k(\boldsymbol{x},t)$.

On the global boundary $\partial \Omega$, boundary conditions are enforced to satisfy Eq.~\eqref{eq:6b} at specific collocation points, yielding:
\begin{equation}
    \boldsymbol{B} \boldsymbol{\beta} = g,
    \label{eq:9}
\end{equation}
where $\boldsymbol{B}$ is the coefficient matrix related to the boundary conditions, and $g$ is the column vector obtained from the discretized boundary conditions in Eq.~\eqref{eq:6b}.

To ensure solution smoothness across shared interfaces $\mathcal{I}_{pq}$, continuity conditions are imposed. 
Substituting the solution expression \eqref{eq:777} into the continuity condition \eqref{eq:5} yields:
\begin{equation}
   \frac{\partial^{\alpha}}{\partial x^{\alpha}_s}\left(\sum_{j=1}^{M}\beta^{p}_{ji}\cdot \phi^{p}_{j} \right) - \frac{\partial^{\alpha}}{\partial x^{\alpha}_s}\left(\sum_{j=1}^{M}\beta^{q}_{ji}\cdot \phi^{q}_{j} \right) =0, \quad \alpha=1,2,\dots,k_s, \quad \forall x \in \mathcal{I}_{pq}, \label{eq:10}
\end{equation}

Enforcing this condition at specific collocation points yields the boundary constraint equation between adjacent subdomains:
\begin{equation}
    \boldsymbol{D}^{pq}_k \boldsymbol{\beta}^p - \boldsymbol{D}^{qp}_k\boldsymbol{\beta}^q = 0,
    \label{eq:11}
\end{equation}
where $\boldsymbol{D}^{pq}_k$ represents the discretized matrix of the basis functions and their derivatives on the shared boundary $\mathcal{I}_{pq}$, ensuring the $C^k$ continuity condition across subdomain interfaces.

Combining Eqs.~\eqref{eq:8}, \eqref{eq:9}, and \eqref{eq:11} yields the global algebraic system:
\begin{equation}
    \begin{bmatrix}
        \boldsymbol{A}^1           &0          &0                           &0           &0\\
        0                          &\ddots     &0                           &0           &0\\
        0                          &0          &\boldsymbol{A}^k            &0           &0\\
        0                          &0          &0                           &\ddots      &0\\
        0                          &0          &0                           &0           &\boldsymbol{A}^{N_k}\\
        \hdashline
        \boldsymbol{B}^1           &0          &0                           &0           &0\\
        0                          &\ddots     &0                           &0           &0\\
        0                          &0          &\boldsymbol{B}^k            &0           &0\\
        0                          &0          &0                           &\ddots      &0\\
        0                          &0          &0                           &0           &\boldsymbol{B}^{N_k}\\
        \hdashline
        \boldsymbol{D}_{0}^{1,2}   &-\boldsymbol{D}_{0}^{2,1}&0             &0           &0\\
        0                          &\ddots     &\ddots                      &0           &0\\
        0                          &0          &\boldsymbol{D}_{0}^{p,q}    &-\boldsymbol{D}_{0}^{q,p}&0\\
        0                          &0          &0                           &\ddots      &\ddots\\
        0                          &0          &0                           &\boldsymbol{D}_{0}^{N_{k-1},N_k}&-\boldsymbol{D}_{0}^{N_k,N_{k-1}}\\
        \hdashline
        \vdots& \vdots&\vdots&\ddots&\vdots\\
        \hdashline
        \boldsymbol{D}_{k}^{1,2}             &-\boldsymbol{D}_{k}^{2,1}&0                       &0           &0\\
        0                          &\ddots     &\ddots                      &0           &0\\
        0                          &0          &\boldsymbol{D}_{k}^{p,q}              &-\boldsymbol{D}_{k}^{q,p}&0\\
        0                          &0          &0                           &\ddots&\ddots\\
        0                          &0          &0                           &\boldsymbol{D}_{k}^{N_{k-1},N_k}&-\boldsymbol{D}_{k}^{N_k,N_{k-1}}\\

    \end{bmatrix}
    \begin{bmatrix}
        \beta^{1}\\
        \vdots\\
        \beta^{k}\\
        \vdots\\
        \beta^{N_k}\\
    \end{bmatrix}
    =
    \begin{bmatrix}
        f^{1}\\
        \vdots\\
        f^{k}\\
        \vdots\\
        f^{N_k}\\
        \hdashline
        g^{1}\\
        \vdots\\
        g^{k}\\
        \vdots\\
        g^{N_k}\\
        \hdashline
        \mathbf{0}\\
        \hdashline
        \vdots\\
        \hdashline
        \mathbf{0}\\
    \end{bmatrix}
    \label{eq:12}
\end{equation}
where the number of rows in $\boldsymbol{A}^k$ is determined by the number of collocation points within the subdomain, and the number of columns equals the subspace dimension $M$. 
The number of rows in $\boldsymbol{B}^k$ is determined by the number of collocation points, and its number of columns is also $M$. 
$\beta^k$ is an $M$-dimensional column vector, and $\mathbf{0}$ in the right-hand side column vector represents a zero vector.

We rewrite Eq.~\eqref{eq:12} as:
\begin{equation}
    \boldsymbol{A}\boldsymbol{\beta} = \boldsymbol{b}.
    \label{eq:13}
\end{equation}
To obtain the global solution, the combination coefficients $\boldsymbol{\beta}$ are determined by solving Eq.~\eqref{eq:13} using the least squares method:
\begin{equation}
    \boldsymbol{\beta} = \mathop{\textit{argmin}}\limits_{\beta} ||\boldsymbol{A}\boldsymbol{\beta} - \boldsymbol{b}||_2^2.
\end{equation}

By unifying the local solutions within subdomains, continuity conditions across subdomain interfaces, boundary conditions, and the source term, the combination coefficients $\boldsymbol{\beta}$ are solved through the least squares method. 
This ensures high accuracy of the solution across the entire domain. The algorithm for solving linear problems is summarized as follows:

\begin{algorithm}[htbp]
    \caption{DD-SNN Algorithm for linear problem}
    \label{alg2:DD-SNN}
    \begin{algorithmic}[1]
        \STATE \textbf{Input}: Domain $\Omega$; number of subdomains $N_k$; maximum number of epochs $N_{max}$; tolerance $\varepsilon$; subspace dimension $M$; network structure $\{d, n^{(1)}, n^{(2)}, \dots, n^{(L+1)}, s\}$.
        
        \STATE \textbf{Step 1}: Partition the domain $\Omega$ into $N_k$ independent, non-overlapping subdomains $\{\Omega_k\}_{k=1}^{N_k}$.
        
        \STATE \textbf{Step 2}: \FOR{each subdomain $\Omega_k$}
        \STATE Initialize parameters $\theta^k$ (weights and biases).
        \STATE Set initial combination coefficients $\boldsymbol{\beta}^k$.
        \STATE Construct the loss function $\mathcal{L}_{PDE}^k(\boldsymbol{x}, t;\theta^k, \boldsymbol{\beta^k})$ and optimize $\theta^k$ by minimizing the loss function.
        \STATE Sample training points $\{(\boldsymbol{x}^i, t^i)\}_{i=1}^{N_r}$ in $\Omega_k$.
        \WHILE{stopping condition not met and the number of epochs $< N_{max}$}
        \STATE Minimize $\mathcal{L}_{PDE}^{k}(\boldsymbol{x}, t; \theta^k, \boldsymbol{\beta}^k)$ to update $\theta^k$.
        \STATE Check if $\frac{\mathcal{L}_{PDE}^{k}}{\mathcal{L}_{PDE_0}^{k}} \leq \varepsilon$.
        \IF{condition is satisfied}
        \STATE Stop training.
        \ENDIF
        \ENDWHILE
        \ENDFOR

        \STATE \textbf{Step 3}: Obtain the basis functions in the subspace layer and represent the solution as a linear combination of these basis functions.

        \STATE \textbf{Step 4}: Enforce the equation at interior points, boundary conditions at boundary points, and $C^k$ continuity conditions on shared interfaces $\mathcal{I}_{pq}$ to construct a global linear algebraic system.
        
        \STATE \textbf{Step 5}: Solve the linear algebraic system to obtain the optimal combination coefficients $\boldsymbol{\beta}$.
        
        \STATE \textbf{Step 6}: Assemble the local solutions into the global solution $\hat{u}(\boldsymbol{x},t)$.
        
        \STATE \textbf{Output}: Numerical solution $\hat{u}(\boldsymbol{x},t)$.
    \end{algorithmic}
\end{algorithm}

\section{Nonlinear Equations}
This section describes the DD-SNN method for solving nonlinear problems. 
We adopt traditional nonlinear iterative methods to address nonlinear problems, specifically using Picard and Newton iteration methods for solving the nonlinear equations \eqref{eq:1}.

\subsection{Picard Iteration}

The core idea of the Picard iteration method is to linearize the nonlinear equation by fixing the nonlinear term $\mathcal{N}(\hat{u}, \hat{u_t},\nabla \hat{u},\nabla^2 \hat{u}, \dots)$. Specifically, an initial guess is chosen, and nonlinear iteration begins. 
After training, an initial solution is naturally obtained and used as the initial guess $\hat{u}^{(0)}$. 
During the nonlinear iteration, the basis functions remain unchanged, and only the combination coefficients $\boldsymbol{\beta}$ are updated to produce the refined solution $\hat{u}$. 
In the $(n+1)$-th iteration, the nonlinear term $\mathcal{N}$ is fixed to the known value from the previous iteration, $\mathcal{N}(\hat{u}^{(n)})$, and substituted into Eq.~\eqref{eq:1a}, resulting in the linearized form:
\begin{equation}
    \mathcal{L} \hat{u}^{(n+1)}(\boldsymbol{x}, t) = f(\boldsymbol{x}, t) - \mathcal{N}(\hat{u}^{(n)})(\boldsymbol{x},t), \quad (\boldsymbol{x}, t) \in \Omega.
    \label{eq:15}
\end{equation}

The solution procedure is as follows. 
To solve Eq.~\eqref{eq:1a}, we follow the steps in Algorithm~\ref{alg2:DD-SNN}. 
For each subdomain $\Omega_k$, the loss function \eqref{eq:3} is minimized to optimize the parameters $\theta^k$, yielding a set of basis functions $\phi_j^k (j=1,2, \dots, M)$. 
These basis functions remain fixed during the subsequent nonlinear iterations, and the solution is represented as a linear combination of them. 
Substituting this solution into the nonlinear equation gives:
\begin{equation}
    \mathcal{L}\left(\sum^{M}_{j=1}\boldsymbol{\beta}^k_{ji}\cdot\phi_j^k(\boldsymbol{x},t) \right)
+\mathcal{N}(\cdot)=f^k(\boldsymbol{x},t), \quad 1\leq i \leq s, \quad (\boldsymbol{x},t) \in \Omega_k,
\label{eq:16}
\end{equation}
where $\mathcal{N}(\cdot)=\mathcal{N}\left(\sum_{j=1}^{M}\boldsymbol{\beta}^k_{ji} \cdot \phi_{j}^k ,\frac{\partial }{\partial t}\left(\sum_{j=1}^{M}\boldsymbol{\beta}^k_{ji} \cdot \phi_{j}^k\right), \frac{\partial}{\partial x}\left(\sum_{j=1}^{M}\boldsymbol{\beta}^k_{ji} \cdot \phi_{j}^k\right),\dots\right)$.

For the nonlinear equation \eqref{eq:16}, we use the Picard iteration method. 
Following the approach for solving linear problems, the following matrix form is obtained in  the subdomain $\Omega_k$,
\begin{equation}
    \boldsymbol{A}^k \boldsymbol{\beta}^{k,(n+1)} = \overline{f}^{k,(n)}(\boldsymbol{\beta}^{k,(n)}),
    \label{eq:17}
\end{equation}
where $\boldsymbol{\beta}^{k,(n+1)}$ is the combination coefficient vector in the $(n+1)$-th iteration for subdomain $\Omega_k$, and $\overline{f}^{k,(n)}$ is a vector dependent on $\boldsymbol{\beta}^{k,(n)}$.

The global boundary conditions are incorporated as:
\begin{equation}
    \boldsymbol{B} \boldsymbol{\beta}^{(n+1)} = g,
\end{equation}
and the continuity conditions are enforced as:
\begin{equation}
    \boldsymbol{D}^{pq}_k \boldsymbol{\beta}^{p,(n+1)} - \boldsymbol{D}^{qp}_k\boldsymbol{\beta}^{q,(n+1)} = 0.
\end{equation}
These conditions yield the global linear algebraic system:
\begin{equation}
    \boldsymbol{A}\boldsymbol{\beta}^{(n+1)} = \boldsymbol{b}(\boldsymbol{\beta}^{(n)}).
    \label{eq:19}
\end{equation}

During the nonlinear iteration, the basis functions remain fixed, and we solve the linear system \eqref{eq:19} using the least squares method to update the combination coefficients $\boldsymbol{\beta}^{(n+1)}$, thereby obtaining $\hat{u}^{(n+1)}$. 
The iteration continues until the convergence condition is satisfied:
\begin{equation}
    ||\hat{u}^{(n+1)} - \hat{u}^{(n)}||_{\infty} \leq \varepsilon_{non},
    \label{eq:20}
\end{equation}
where $\varepsilon_{non}$ is a predefined convergence threshold.

The Picard iteration method is simple to implement and effectively solves nonlinear equations. However, when the nonlinear term is strong, the method may face convergence issues.

\subsection{Newton Iteration}

To enhance the convergence of iterations, we introduce the Newton iteration method. 
Similar to the Picard method, an initial guess $\boldsymbol{\beta}^{(0)}$ is selected, yielding the initial solution $\hat{u}^{(0)}(\boldsymbol{x},t)$. Newton iterations are sensitive to the initial guess, requiring an appropriate initial value to ensure convergence. 
In the DD-SNN method, since the training for each subdomain is entirely independent, the solution obtained after training might not yet approximate the solution of equation; 
boundary conditions and inter-subdomain continuity must be applied. 
Thus, the Picard iteration method is first employed to obtain the initial guess for Newton iteration.

For the nonlinear equation \eqref{eq:1a}, the nonlinear residual is defined as:
\begin{equation}
    R(\hat{u}) = \mathcal{L}\hat{u}(\boldsymbol{x}, t) + \mathcal{N}(\hat{u}, \hat{u_t},\nabla \hat{u},\nabla^2 \hat{u}, \dots) - f(\boldsymbol{x}, t),
    \label{eq:21}
\end{equation}
where the nonlinear equation \eqref{eq:1a} is equivalent to $R(\hat{u})=0$. Substituting the expression for $\hat{u}$ into the above equation gives:
\begin{equation}
    R\left(\sum^{M}_{j=1}\boldsymbol{\beta}_{ji}\cdot\phi_j(\boldsymbol{x},t)\right) = \mathcal{L}\left(\sum^{M}_{j=1}\boldsymbol{\beta}_{ji}\cdot\phi_j(\boldsymbol{x},t)\right) + \mathcal{N}(\cdot) - f(\boldsymbol{x}, t),
\end{equation}
where $\mathcal{N}(\cdot) = \mathcal{N}\left(\sum_{j=1}^{M}\boldsymbol{\beta}_{ji} \cdot \phi_{j} ,\frac{\partial }{\partial t}\left(\sum_{j=1}^{M}\boldsymbol{\beta}_{ji} \cdot \phi_{j}\right), \frac{\partial}{\partial x}\left(\sum_{j=1}^{M}\boldsymbol{\beta}_{ji} \cdot \phi_{j}\right),\dots\right)$.

The original problem is thus transformed into a nonlinear problem with respect to $\boldsymbol{\beta}$. 
In the $(n+1)$-th iteration, the nonlinear residual \eqref{eq:21} is linearized using a Taylor expansion around the solution $\boldsymbol{\beta}^{(n)}$:
\begin{equation}
    R(\boldsymbol{\beta}^{(n+1)}) \approx R(\boldsymbol{\beta}^{(n)}) + J(\boldsymbol{\beta}^{(n)})\Delta \boldsymbol{\beta},
    \label{eq:22}
\end{equation}
where $J(\boldsymbol{\beta}^{(n)}) =\frac{\partial}{\partial \boldsymbol{\beta}^{(n)}}\left(\mathcal{L}(\boldsymbol{\beta}^{(n)}\cdot\phi)+\mathcal{N}(\boldsymbol{\beta}^{(n)}\cdot\phi)\right)$ is the Jacobian matrix, and $\Delta \boldsymbol{\beta} = \boldsymbol{\beta}^{(n+1)} - \boldsymbol{\beta}^{(n)}$ is the correction. Setting $R(\boldsymbol{\beta}^{(n+1)})=0$ in the expansion yields the linearized equation:
\begin{equation}
    J(\boldsymbol{\beta}^{(n)})\Delta \boldsymbol{\beta} = -R(\boldsymbol{\beta}^{(n)}).
    \label{eq:23}
\end{equation}

To solve the linearized equation \eqref{eq:23}, we follow the steps in Algorithm~\ref{alg2:DD-SNN}. 
For each subdomain $\Omega_k$, the equation becomes:
\begin{equation}
    J^k\left(\sum^M_{j=1}\boldsymbol{\beta}_{ji}^{k,(n)}\cdot\phi_j^k(\boldsymbol{x},t)\right)\Delta \boldsymbol{\beta}^{k} = -R^k\left(\sum^M_{j=1}\boldsymbol{\beta}_{ji}^{k,(n)}\cdot\phi_j^k(\boldsymbol{x},t)\right),
    \label{eq:24}
\end{equation}
where $J^k$ is the Jacobian matrix on subdomain $\Omega_k$, $\Delta \boldsymbol{\beta}^k=\boldsymbol{\beta}^{k,(n+1)}-\boldsymbol{\beta}^{k,(n)}$ is the correction for $\boldsymbol{\beta}$ on $\Omega_k$, and $R^k$ is the residual vector on $\Omega_k$.

Enforcing \eqref{eq:24} at collocation points results in the matrix form:
\begin{equation}
    \boldsymbol{A}^k(\boldsymbol{\beta}^{k,(n)}) \Delta \boldsymbol{\beta}^k = -\boldsymbol{r}^k(\boldsymbol{\beta}^{k,(n)}),
    \label{eq:25}
\end{equation}
where $\boldsymbol{A}^k$ is the coefficient matrix on subdomain $\Omega_k$, and $\boldsymbol{r}^k$ is the discretized vector of the residual $R^k$.

Since the Newton method iterates over the correction $\Delta \boldsymbol{\beta}$, the boundary and continuity conditions must also be expressed in terms of $\Delta \boldsymbol{\beta}$, yielding:
\begin{equation}
    \boldsymbol{B} \Delta \boldsymbol{\beta} = g - \boldsymbol{B} \cdot\boldsymbol{\beta}^{(n)},
    \label{eq:26}
\end{equation}
\begin{equation}
    (\boldsymbol{D}^{pq}_k - \boldsymbol{D}^{qp}_k)\Delta \boldsymbol{\beta} = -(\boldsymbol{D}^{pq}_k - \boldsymbol{D}^{qp}_k)\boldsymbol{\beta}^{(n)},
    \label{eq:27}
\end{equation}
where $\boldsymbol{\beta}^{(n)}$ is the global coefficient vector at the $n$-th iteration, and $\Delta \boldsymbol{\beta} = \boldsymbol{\beta}^{(n+1)}-\boldsymbol{\beta}^{(n)}$.

Combining Eqs.~\eqref{eq:25}, \eqref{eq:26}, and \eqref{eq:27} yields the global algebraic system:
\begin{equation}
    \boldsymbol{A}(\boldsymbol{\beta}^{(n)})\Delta \boldsymbol{\beta} = \boldsymbol{b}(\boldsymbol{\beta}^{(n)}).
    \label{eq:28}
\end{equation}
The correction $\Delta \boldsymbol{\beta}$ is solved using the linear least squares method, and the coefficients are updated as:
\begin{equation}
    \boldsymbol{\beta}^{(n+1)} = \boldsymbol{\beta}^{(n)} + \Delta \boldsymbol{\beta}.
\end{equation}
The new iterative solution $\hat{u}^{(n+1)}$ is obtained. The iteration stops once the convergence condition for \eqref{eq:16} is satisfied.

The Newton iteration method involves the computation of the Jacobian matrix, resulting in higher computational cost per iteration compared to the Picard method. However, Newton iteration generally achieves second-order convergence, allowing for faster convergence, especially for strongly nonlinear problems, offsetting its higher per-iteration cost.

\section{Numerical Experiments}

This section presents numerical experiments to validate the performance of the proposed DD-SNN method. 
These experiments include steady-state and time-dependent linear and nonlinear differential equations. 
Unless otherwise specified, all experiments use the Tanh function as the activation function, the Adam optimizer with a learning rate of 0.001, and a global random seed of 202 to ensure reproducibility. 
For the ELM and LocELM methods, initialization parameters are uniformly distributed in the interval $[-1,1]$. Except for the ELM method, which does not include hidden layers, all methods use the same network structure.

In the DD-SNN method, computational cost refers to the total time of solving PDE, including the training of basis functions, automatic differentiation (e.g., $\Phi_j^{k},\frac{\partial}{\partial \boldsymbol{x}} \Phi_j^{k}$), the time to construct the coefficient matrix and least-squares problem, and the time for nonlinear iteration. 
Epochs represent the average number of training per subdomain.

\subsection{Linear Example}
\subsubsection{One-Dimensional Linear Helmholtz Equation}

We first consider the boundary value problem for a one-dimensional linear Helmholtz equation on the domain $\Omega = [a,b]$:
\begin{subequations}
    \begin{align}
        \frac{\mathrm{d}^2 u}{\mathrm{d}x^2} - \lambda u &= f(x), \label{eq:33a} \\
        u(a) &= h_1, \label{eq:33b} \\
        u(b) &= h_2, \label{eq:33c}
    \end{align}
    \label{eq:33}
\end{subequations}
where $u(x)$ is the solution to be determined, $f(x)$ is the specified source term, and $h_1$ and $h_2$ are the given boundary values. 
In this experiment, we set $\lambda = 10$, $a = 0$, and $b = 8$. The source term $f(x)$ is chosen such that the exact solution of Eq.~\eqref{eq:33} is:
\begin{equation}
    u(x) = \sin(3 \pi x + \frac{3\pi}{20}) \cos(2\pi x+\frac{\pi}{10})+2.
\end{equation}
The boundary conditions $h_1$ and $h_2$ are determined from the exact solution.

To ensure fairness, all methods are tested on 800 uniformly sampled points. 
In this experiment, the PINNs algorithm is implemented using the DEEPXDE package, with a fully connected neural network containing three hidden layers, each with 100 neurons, and trained for 50,000 epochs. 
The DGM and DRM methods use the same network hyperparameters for consistency. 
The ELM method uses a single-layer network with 100 neurons. 
The LocELM method divides the domain into four subdomains, using three hidden layers with 100 neurons each. For the DD-SNN method, the domain is also divided into four subdomains (all numerical experiments in this example are based on this setting unless otherwise specified). 
Each subdomain uses a neural network with two hidden layers of 100 neurons each and a subspace dimension of 100. 
To ensure comparability, all methods adopt a uniform sampling strategy.

Table~\ref{tab:1} presents the errors and epochs of different neural network methods for solving Eqs.~\eqref{eq:33a}-\eqref{eq:33c}. 
The results show that after 50,000 epochs, the PINNs, DGM, and DRM methods still exhibit significant errors, and increasing the number of epochs does not significantly improve accuracy. 
The ELM method yields large errors, while the LocELM method achieves an $L^2$ error of $4.96 \times 10^{-9}$. 
The DD-SNN method achieves an accuracy of $4.98 \times 10^{-11}$ after only 630 epochs, demonstrating superior accuracy compared to other methods.

\begin{table}[htbp]
    \setlength{\abovecaptionskip}{0cm}
    \setlength{\belowcaptionskip}{0.2cm}
    \centering
    \caption{Comparison of errors, epochs, and computational efficiency for different methods on the linear Helmholtz equation}
    \begin{tabular}{ccccc}
    \hline
    Method & $\|e\|_{L^2}$ & $\|e\|_{L^{\infty}}$ & Epochs  & CPU time(s) \\ \hline
    PINN   & $1.05e-02$    & $5.63e-02$           & $50000$ &612.20\\
    DGM    & $8.82e-03$    & $2.22e-02$           & $50000$ &561.71\\
    DRM    & $1.29e-02$    & $3.92e-02$           & $50000$ &619.43\\
    ELM    & $3.90e-01$    & $8.69e-01$           & $-    $ &0.02\\
    LocELM & $4.96e-09$    & $1.44e-08$           & $-    $ &0.05\\
    DD-SNN & $4.98e-11$    & $2.02e-10$           & $630$  &18.48\\ \hline
    \end{tabular}
    \label{tab:1}
\end{table}

\begin{figure}
    \centering
    \includegraphics[width=0.6\textwidth]{ 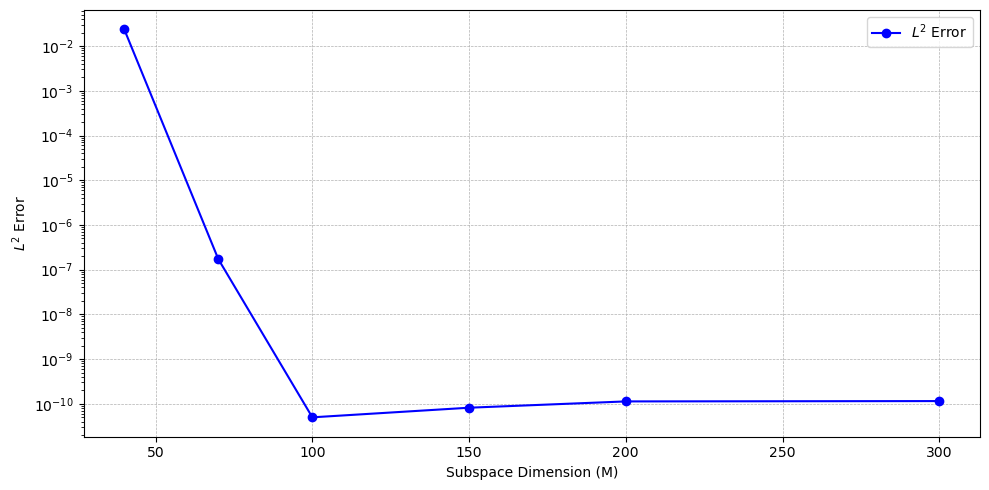}
    \caption{Impact of subspace dimension on  $L^2$ error}
    \label{fig:2}
\end{figure}

Figure~\ref{fig:2} shows the effect of subspace dimensions on errors at 300 uniform points. 
Numerical results indicate that as the number of subspaces increases, the errors decrease gradually. 

Next, we investigate the effect of the number of subdomains on the error. 
Keeping the total number of uniformly sampled points fixed at 800, the domain is divided into varying numbers of subdomains. 
When the number of subdomains is 1, the DD-SNN method degenerates into the standard SNN method. 
Table~\ref{tab:3} presents the impact of subdomain numbers on errors and computational efficiency. 
The results show that as the number of subdomains increases, the error decreases and the number of epochs are reduced. 
This is because with more subdomains, the global solution is decomposed into relatively simpler local solutions, reducing the optimization complexity for each subdomain. 
However, when the number of subdomains becomes excessively large, the number of sampling points within each subdomain decreases, leading to an increase in error.
When the number of subdomains is 1, the number of epochs reach the maximum and the accuracy is not satisfied.  
This is because the subspace's ability to approximate the solution space is insufficient, and further training can improve the accuracy of the method.

\begin{table}[htbp]
    \setlength{\abovecaptionskip}{0cm}
    \setlength{\belowcaptionskip}{0.2cm}
    \centering
    \caption{Impact of subdomain numbers on errors and training efficiency}
    \begin{tabular}{ccccc}
    \hline
    Subdomain    & $\|e\|_{L^2}$ & $\|e\|_{L^{\infty}}$ & Epochs   & CPU time(s)\\ \hline
    1            & $2.99e-01$    & $1.24e+00$           & $5000$   & 72.91\\
    2            & $1.66e-05$    & $7.87e-05$           & $2407$   & 68.29\\
    4            & $4.98e-11$    & $2.02e-10$           & $630$    & 19.07\\
    6            & $3.46e-11$    & $1.51e-10$           & $482$    & 18.23\\
    8            & $2.73e-11$    & $9.39e-11$           & $326$    & 14.67\\
    16           & $9.49e-13$    & $3.13e-12$           & $126$    & 10.18\\ 
    32           & $2.55e-11$    & $9.25e-11$           & $55$     & 8.31\\\hline
    \end{tabular}
    \label{tab:3}
\end{table}

We further investigated the impact of the number of hidden layers on accuracy and computational efficiency. 
Keeping the number of subdomains fixed at 4 and each hidden layer containing 100 neurons, the number of hidden layers was varied. 
Although this configuration may not be optimal, it highlights the advantages of our method. 
Table~\ref{tab:4} shows that the DD-SNN method achieves high accuracy regardless of the number of hidden layers, although there are slight differences in accuracy.

\begin{table}[htbp]
    \setlength{\abovecaptionskip}{0cm}
    \setlength{\belowcaptionskip}{0.2cm}
    \centering
    \caption{Impact of hidden layer numbers on errors, training epochs, and computational efficiency}
    \begin{tabular}{ccccc}
    \hline
    Hidden layer & $\|e\|_{L^2}$ & $\|e\|_{L^{\infty}}$ & Epochs  & CPU time(s)\\ \hline
    0            & $1.70e-07$    & $7.06e-07$           & $2823$  & 122.97\\
    1            & $1.24e-07$    & $2.93e-07$           & $1180$  & 24.87\\
    2            & $4.98e-11$    & $2.02e-10$           & $630 $  & 19.34\\
    3            & $2.98e-10$    & $1.15e-09$           & $1313$  & 52.71\\
    4            & $2.27e-10$    & $1.13e-09$           & $956$   & 48.10\\ 
    5            & $8.38e-10$    & $3.73e-09$           & $709$   & 42.41\\ 
    6            & $4.24e-10$    & $1.17e-09$           & $733$   & 49.99\\
    7            & $2.20e-09$    & $8.44e-09$           & $757$   & 58.94\\
    8            & $2.36e-09$    & $1.20e-08$           & $460$   & 40.58\\
    \hline
    \end{tabular}
    \label{tab:4}
\end{table}

Next, we tested the effect of three different sampling methods on accuracy.  
The results in Table~\ref{tab:5} show that Gaussian sampling and uniform sampling achieve similar accuracy, with slightly lower errors compared to random sampling. 
The difference, however, is within one order of magnitude. 
Gaussian sampling uses Gauss-Lobatto-Legendre points, while random sampling includes subdomain endpoints to satisfy continuity conditions. These results demonstrates that the sampling method has a limited impact on the final results, indicating the robustness of the DD-SNN method.

\begin{table}[htbp]
    \setlength{\abovecaptionskip}{0cm}
    \setlength{\belowcaptionskip}{0.2cm}
    \centering
    \caption{Impact of sampling methods on errors, epochs, and computational efficiency}
    \renewcommand{\arraystretch}{1.3} 
    \begin{tabular}{ccccc}
    \hline
    Sampling method & $\|e\|_{L^2}$ & $\|e\|_{L^{\infty}}$ & Epochs  & CPU time(s)\\  \hline
    Uniform         & $4.98e-11$    & $2.02e-10$           & $630$   &19.07\\
    Gaussian        & $9.43e-11$    & $2.66e-10$           & $647$   &33.91\\
    Random          & $6.99e-09$    & $2.78e-08$           & $1792$  &57.16\\
    \hline
    \end{tabular}
    \label{tab:5}
\end{table}

Furthermore, to avoid experimental randomness, we tested the effects of different random seeds. 
Table~\ref{tab:6} shows the results for various random seeds. 
It can be observed that different random seeds have little impact on the results, demonstrating that the DD-SNN method is insensitive to the choice of random seed and exhibits robust performance.

\begin{table}[htbp]
    \setlength{\abovecaptionskip}{0cm}
    \setlength{\belowcaptionskip}{0.2cm}
    \centering
    \caption{Impact of random seeds on errors, epochs, and computational efficiency}
    \begin{tabular}{ccccc}
    \hline
    Random seed  & $\|e\|_{L^2}$ & $\|e\|_{L^{\infty}}$ & Epochs  & CPU time(s)\\  \hline
    1            & $3.24e-10$    & $8.78e-10$           & $553$  &17.96\\
    202          & $4.98e-11$    & $2.02e-10$           & $630$  &20.79\\
    666          & $7.52e-11$    & $2.31e-10$           & $519$  &17.15\\
    888          & $1.40e-10$    & $5.66e-10$           & $728$  &24.12\\
    2024         & $5.97e-11$    & $2.60e-10$           & $689$  &22.79\\
    666888       & $5.76e-10$    & $2.59e-09$           & $771$  &25.16\\
    \hline
    \end{tabular}
    \label{tab:6}
\end{table}

\begin{figure}
    \centering
    \includegraphics[width=0.8\textwidth]{ 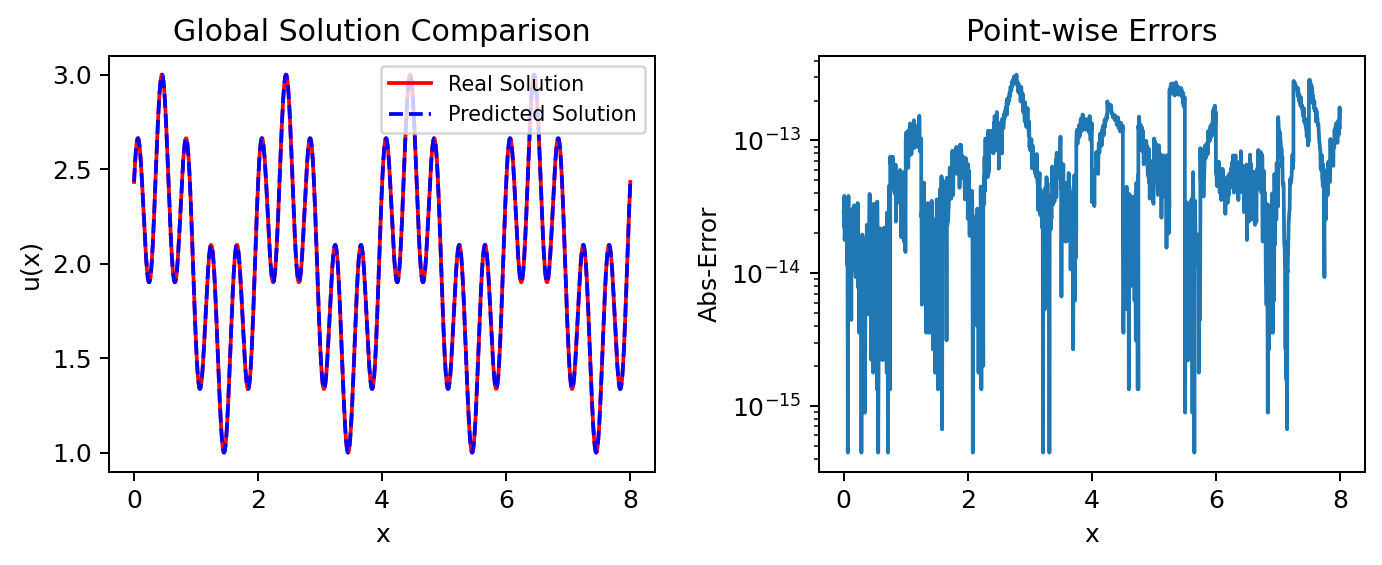}
    \caption{Comparison of numerical and exact solutions for the one-dimensional Helmholtz equation}
    \label{fig:3}
\end{figure}

Figure~\ref{fig:3} illustrates the comparison between the numerical and exact solutions obtained by the DD-SNN method for the one-dimensional Helmholtz equation \eqref{eq:33a}-\eqref{eq:33c}. 
The domain is divided into 32 subdomains, with 50 uniformly sampled points in each subdomain. 
The pointwise error distribution is also displayed.

\subsubsection{Two-Dimensional Linear Poisson Equation}

We consider the two-dimensional Poisson equation on the domain $\Omega = [a,b] \times [a,b]$:
\begin{subequations}
    \begin{alignat}{2}
        -\Delta u(x,y) &= f(x,y), &\quad &(x, y) \in \Omega, \label{eq:35a} \\
        u(x,y)         &= g(x,y), &\quad &(x, y) \in \partial \Omega,   \label{eq:35b}
    \end{alignat}
    \label{eq:35}
\end{subequations}
where $u(x,y)$ is the solution to be determined, $f(x,y)$ is the source term, and $g(x,y)$ represents the boundary conditions. 
We set $a=0$ and $b=2$. The source term $f(x,y)$ is chosen such that the exact solution is:
\begin{equation}
    u(x,y) = \sin(\pi x) \sin(\pi y).
\end{equation}

Table~\ref{tab:7} presents the errors, epochs, and computational efficiency for different methods applied to Eq.\eqref{eq:35}. 
All methods use $30\times30\times16$ uniform sampling points and the results show that after $50,000$ epochs, the PINNs, DGM, and DRM methods still exhibit significant errors, and further increasing the number of epochs does not significantly improve their accuracy. 
The ELM method, with a single layer of $300$ neurons, achieves an $L^2$ error of $6.39 \times 10^{-9}$. The LocELM method fails to converge under the same hyperparameter settings as DD-SNN, resulting in poor accuracy.

The DD-SNN method divides the domain into $16$ subdomains. 
The hyperparameters include two hidden layers per subdomain, each with $100$ neurons, and a subspace dimension of $300$. 
Compared to other methods, the DD-SNN method significantly reduces the total number of epochs and demonstrates superior accuracy.

\begin{table}[htbp]
    \setlength{\abovecaptionskip}{0cm}
    \setlength{\belowcaptionskip}{0.2cm}
    \centering
    \caption{Comparison of errors, training epochs, and computational efficiency for different methods on the two-dimensional Poisson equation}
    \begin{tabular}{ccccc}
    \hline
    Method & $\|e\|_{L^2}$ & $\|e\|_{L^{\infty}}$ & Epochs  & CPU time (s)\\ \hline
    PINN   & $1.25e-02$    & $2.27e-02$           & $50000$ & 16146.89\\
    DGM    & $2.58e-03$    & $6.72e-03$           & $50000$ & 22534.92\\
    DRM    & $8.36e-03$    & $1.34e-02$           & $50000$ & 22644.70\\
    ELM    & $6.39e-09$    & $5.31e-08$           & $-    $ & 14.61\\
    LocELM & $4.15e-01$    & $1.29e+00$           & $-    $ & 60.60\\
    DD-SNN & $4.94e-12$    & $6.61e-11$           & $57$    & 154.61\\ \hline
    \end{tabular}
    \label{tab:7}
\end{table}

\begin{figure}
    \centering
    \includegraphics[width=0.6\textwidth]{ 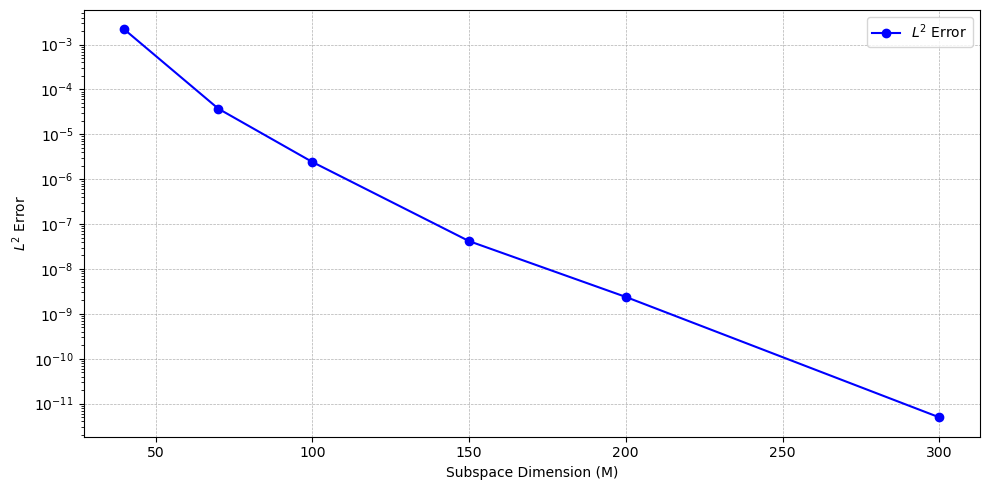}
    \caption{Impact of subspace dimension on $L^2$ error}
    \label{fig:4}
\end{figure}

Figure~\ref{fig:4} shows the impact of subspace dimensions on the results at $30\times30\times16$ uniform points. 
Numerical results indicate that as the number of subspaces increases, the error gradually decreases.

Table~\ref{tab:9} shows that as the number of subdomains increases, the error gradually decreases, indicating that finer partitions can effectively improve solution accuracy. 
When the number of subdomains increases from 1 to 16, the $L^2$ error decreases from $7.92 \times 10^{-7}$ to $4.94 \times 10^{-12}$, while the $L^{\infty}$ error decreases from $3.40 \times 10^{-6}$ to $6.61 \times 10^{-11}$. 
This demonstrates that increasing the number of subdomains significantly enhances solution accuracy.

\begin{table}[htbp]
    \setlength{\abovecaptionskip}{0cm}
    \setlength{\belowcaptionskip}{0.2cm}
    \centering
    \caption{Impact of subdomain numbers on errors and training efficiency}
    \begin{tabular}{ccccc}
    \hline
    Subdomain    & $\|e\|_{L^2}$ & $\|e\|_{L^{\infty}}$ & Epochs   & CPU time(s)\\ \hline
    1            & $7.92e-07$    & $3.40e-06$           & $774$    & 33.06\\
    4            & $1.46e-10$    & $2.27e-09$           & $288$    & 189.01\\
    9            & $2.40e-11$    & $5.38e-10$           & $90$     & 115.95\\
    16           & $4.94e-12$    & $6.61e-11$           & $57$     & 154.61\\ 
    25           & $4.20e-12$    & $4.45e-11$           & $50$     & 296.70\\ \hline
    \end{tabular}
    \label{tab:9}
\end{table}

Next, we analyze the impact of the number of hidden layers on errors and training costs. 
The global sampling points are set to $30 \times 30 \times 16$, and other network hyperparameters are fixed, including 100 neurons per hidden layer, 16 subdomains, and a subspace dimension of 300.

As shown in Table~\ref{tab:10}, the number of hidden layers has a slight impact on training costs, but all configurations converge to good accuracy. 
When the number of hidden layers is 0, the $L^2$ and $L^{\infty}$ errors are $2.74 \times 10^{-13}$ and $4.06 \times 10^{-12}$, respectively, with 271 epochs and a computation time of 186.58 seconds. 
As the number of hidden layers increases to 2, the number of epochs and computation time decrease, reaching a minimum at 2 hidden layers. 
However, the error increases slightly, with $L^2$ and $L^{\infty}$ errors rising to $4.94 \times 10^{-12}$ and $6.61 \times 10^{-11}$, respectively. 
Further increases in hidden layers result in higher errors and training costs. Overall, 2 hidden layers achieve a good balance between errors and training costs.

\begin{table}[htbp]
    \setlength{\abovecaptionskip}{0cm}
    \setlength{\belowcaptionskip}{0.2cm}
    \centering
    \caption{Impact of hidden layer numbers on errors, epochs, and computational efficiency}
    \begin{tabular}{ccccc}
    \hline
    Hidden layer & $\|e\|_{L^2}$ & $\|e\|_{L^{\infty}}$ & Epochs & CPU time(s)\\ 
    \hline
    0            & $2.74e-13$    & $4.06e-12$           & $271$ & 186.58\\
    1            & $1.89e-12$    & $2.56e-11$           & $67$  & 222.78\\
    2            & $4.94e-12$    & $6.61e-11$           & $57$  & 154.61\\
    3            & $1.90e-11$    & $2.30e-10$           & $59$  & 299.63\\
    4            & $3.12e-11$    & $4.84e-10$           & $91$  & 540.04\\ 
    5            & $9.61e-11$    & $2.12e-09$           & $87$  & 541.96\\ 
    6            & $1.20e-09$    & $1.75e-08$           & $76$  & 312.49\\
    7            & $4.61e-10$    & $6.07e-09$           & $74$  & 410.15\\
    8            & $8.38e-09$    & $1.14e-07$           & $87$  & 805.95\\
    \hline
    \end{tabular}
    \label{tab:10}
\end{table}

We also tested the impact of three different sampling methods on accuracy. 
Numerical experiments were conducted under fixed network hyperparameters (100 neurons per hidden layer, 16 subdomains, 2 hidden layers, and a subspace dimension of 300).

As shown in Table~\ref{tab:11}, Gaussian and uniform sampling achieve similar accuracy, both showing good precision, while random sampling yields errors one order of magnitude higher. 
However, all three methods achieve good results for this problem.

\begin{table}[htbp]
    \setlength{\abovecaptionskip}{0cm}
    \setlength{\belowcaptionskip}{0.2cm}
    \centering
    \caption{Impact of sampling methods on errors, epochs, and computational efficiency}
    \renewcommand{\arraystretch}{1.3} 
    \begin{tabular}{ccccc}
    \hline
    Sampling method  & $\|e\|_{L^2}$ & $\|e\|_{L^{\infty}}$ & Epochs  & CPU time(s)\\  \hline
    Uniform          & $4.94e-12$    & $2.02e-10$           & $57$    & 152.22\\
    Gaussian         & $7.20e-12$    & $9.45e-11$           & $51$    & 161.06\\
    Random           & $4.29e-11$    & $2.22e-10$           & $58$    & 136.95\\
    \hline
    \end{tabular}
    \label{tab:11}
\end{table}

We further tested the impact of different random seeds on the results. 
Table~\ref{tab:12} shows that variations in random seeds have minimal effect on accuracy, demonstrating that the DD-SNN method is insensitive to the choice of random seed and exhibits robust performance.

\begin{table}[htbp]
    \setlength{\abovecaptionskip}{0cm}
    \setlength{\belowcaptionskip}{0.2cm}
    \centering
    \caption{Impact of random seeds on errors, epochs, and computational efficiency}
    \begin{tabular}{ccccc}
    \hline
    Random seed  & $\|e\|_{L^2}$ & $\|e\|_{L^{\infty}}$ & Epochs        & CPU time(s)\\  \hline
    1            & $3.20e-12$    & $4.87e-11$           & $57$      &222.69\\
    202          & $4.94e-12$    & $2.02e-10$           & $57$      &154.61\\
    666          & $4.74e-12$    & $7.41e-11$           & $72$      &345.50\\
    888          & $3.64e-12$    & $5.66e-10$           & $61$      &252.33\\
    2024         & $4.09e-12$    & $5.33e-11$           & $72$      &401.58\\
    666888       & $4.82e-12$    & $9.83e-11$           & $68$      &316.11\\
    \hline
    \end{tabular}
    \label{tab:12}
\end{table}

\begin{figure}
    \centering
    \includegraphics[width=0.8\textwidth]{ 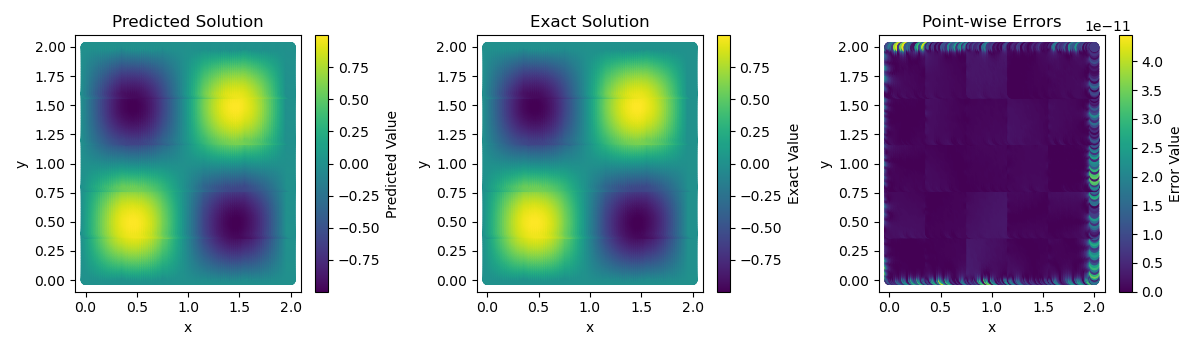} 
    \caption{Comparison of numerical and exact solutions for the two-dimensional Poisson equation} 
    \label{fig:5} 
\end{figure}

Figure~\ref{fig:5} illustrates the numerical and exact solutions of the two-dimensional Poisson equation solved by DD-SNN when the number of subdomains is 25 (refer to Table~\ref{tab:9}). 
The pointwise error distribution is also shown.

\subsubsection{One-Dimensional Parabolic Equation}

We consider the parabolic equation on the spatiotemporal domain $\Omega = [a,b] \times [0,T]$:
\begin{subequations}
    \begin{align}
        \frac{\partial u}{\partial t} - \frac{\partial^2 u}{\partial x^2} &= f(x,t), \label{eq:37a} \\
        u(a,t) &= g_1(t), \label{eq:37b} \\
        u(b,t) &= g_2(t), \label{eq:37c} \\
        u(x,0) &= h(x), \label{eq:37d}
    \end{align}
\end{subequations}
where $u(x,t)$ is the solution to be determined, $f(x,t)$ is the source term, and $g_1(t)$, $g_2(t)$, and $h(x)$ are given functions. 
The source term $f(x,t)$ is chosen such that the exact solution is:
\begin{equation} 
    u(x,t) = 2 e^{-t} \sin(\pi x). 
\end{equation}
In this numerical experiment, we set $a=0$, $b=2$, and $T=2$.

Table~\ref{tab:13} presents the errors and epochs for different neural network methods when solving the one-dimensional parabolic equation. 
The results show that PINNs, DGM, and DRM methods fail to achieve high accuracy even after 50,000 epochs. 
For instance, the $L^2$ error of PINNs is $9.96 \times 10^{-3}$, while DGM and DRM achieve errors of $8.49 \times 10^{-3}$ and $7.82 \times 10^{-3}$, respectively, indicating that these methods struggle to achieve high accuracy even with increased epochs.

The ELM method, using $14,400$ uniform sampling points and a network with $300$ neurons, achieves an $L^2$ error of $1.78 \times 10^{-9}$. 
However, the LocELM method fails to converge under the same hyperparameter settings as DD-SNN, resulting in poor accuracy.

In contrast, the DD-SNN method divides the solution domain into $9$ subdomains and uses a network configuration with two hidden layers, each containing $100$ neurons, and a subspace dimension of $300$. 
With an average of 135 epochs, the DD-SNN method achieves an $L^2$ error of $2.33 \times 10^{-11}$, significantly outperforming PINNs, DGM, and DRM methods. 
This demonstrates that the DD-SNN method, through domain decomposition, achieves significantly higher accuracy and requires fewer epochs to obtain a high-accuracy solution.

\begin{table}[htbp]
    \setlength{\abovecaptionskip}{0cm}
    \setlength{\belowcaptionskip}{0.2cm}
    \centering
    \caption{Comparison of errors, epochs, and computational efficiency for different methods solving the one-dimensional parabolic equation}
    \begin{tabular}{ccccc}
    \hline
    Method & $\|e\|_{L^2}$ & $\|e\|_{L^{\infty}}$ & Epochs         &CPU time(s)\\ \hline
    PINN   & $2.29e-02$    & $7.04e-02$           & $50000$        & 5904.13\\
    DGM    & $8.49e-03$    & $1.15e-02$           & $50000$        & 10200.71\\
    DRM    & $7.82e-03$    & $1.04e-02$           & $50000$        & 1690.05\\
    ELM    & $1.78e-09$    & $1.46e-08$           & $-    $        & 16.80\\
    LocELM & $4.51e-01$    & $2.05e+00$           & $-    $        & 44.82\\
    DD-SNN & $2.33e-11$    & $1.68e-10$           & $98$           & 66.10\\ \hline
    \end{tabular}  
    \label{tab:13}
\end{table}

\begin{figure}
    \centering
    \includegraphics[width=0.6\textwidth]{ 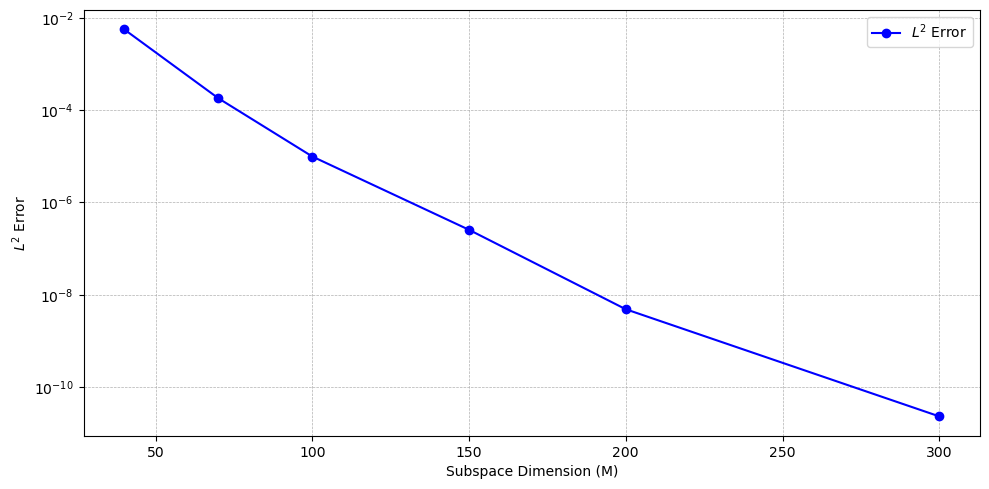}
    \caption{Impact of subspace dimension on $L^2$ error}
    \label{fig:6}
\end{figure}

Figure~\ref{fig:6} shows the effect of different subspace dimensions on errors and computational efficient at $30\times30\times9$ uniform points. 
As the subspace dimension $M$ increases, the $L^2$ error decreases steadily. 
This trend indicates that higher subspace dimensions improve solution accuracy. 
Moreover, the number of epochs remain relatively stable across most configurations, ranging from 73 to 129, demonstrating the computational stability of DD-SNN.

Next, we tested the impact of different numbers of subdomains on solution accuracy. 
By dividing the solution domain into varying numbers of subdomains, Table~\ref{tab:15} presents the effect of subdomain numbers on errors and epochs. 
As the number of subdomains increases, the $L^2$ error gradually decreases, and the number of epochs also reduces. 
For example, when the number of subdomains increases from 1 to 25, the $L^2$ error decreases from $5.18 \times 10^{-10}$ to $1.46 \times 10^{-11}$, and the epochs drop from 364 to 64. 
Increasing the number of subdomains simplifies the global solution into smaller, more manageable local solutions, reducing the training complexity for each subdomain. 
With 9 subdomains, the error is already relatively small, and the training time is relatively low. 
To balance computational cost and accuracy, we use 9 subdomains in subsequent experiments, maintaining a subspace dimension of 300 and global sampling points of $30 \times 30 \times 9$.

\begin{table}[htbp]
    \setlength{\abovecaptionskip}{0cm}
    \setlength{\belowcaptionskip}{0.2cm}
    \centering
    \caption{Impact of subdomain numbers on errors and computational efficiency}
    \begin{tabular}{ccccc}
    \hline
    Subdomain    & $\|e\|_{L^2}$ & $\|e\|_{L^{\infty}}$ & Epochs & CPU time(s)\\ \hline
    1            & $5.18e-10$    & $2.10e-09$           & $364$  & 170.33\\
    4            & $6.97e-11$    & $1.06e-09$           & $108$  & 48.56\\
    9            & $2.33e-11$    & $1.68e-10$           & $98$   & 66.10\\
    16           & $1.62e-11$    & $1.93e-10$           & $72$   & 76.77\\ 
    25           & $1.46e-11$    & $1.14e-10$           & $64$   & 143.49\\ \hline
    \end{tabular} 
    \label{tab:15}
\end{table}

Next, we explored the effect of different numbers of hidden layers on solution accuracy and computational efficiency. 
In this experiment, all configurations used two hidden layers with 100 neurons per layer, 9 subdomains, a subspace dimension of 300, and global uniform sampling points of $30 \times 30 \times 9$, is used. 
Although this configuration does not yield the best results in testing, it highlights the advantages of DD-SNN under different configurations.

Table~\ref{tab:16} shows the impact of hidden layer numbers on $L^2$ error, $L^\infty$ error, epochs, and computational efficiency. 
As the number of hidden layers increases, the error and epochs exhibit some fluctuations, but they all converge to relatively high accuracy while maintaining training time within a reasonable range. 
Although increasing the number of hidden layers enhances the network's expressive capacity, the solution in this case is relatively simple and does not require a very deep network structure. 
Based on the analysis above, we continue to use a two-hidden-layer configuration in subsequent experiments to achieve good accuracy while maintaining low computational costs.

\begin{table}[htbp]
    \setlength{\abovecaptionskip}{0cm}
    \setlength{\belowcaptionskip}{0.2cm}
    \centering
    \caption{Impact of hidden layer numbers on errors, epochs, and computational efficiency}
    \begin{tabular}{ccccc}
    \hline
    Hidden layer & $\|e\|_{L^2}$ & $\|e\|_{L^{\infty}}$ & Epochs         & CPU time(s)\\ \hline
    0            & $6.00e-12$    & $4.74e-11$           & $242$   & 40.30\\
    1            & $1.55e-11$    & $1.51e-10$           & $77$    & 37.88\\
    2            & $2.33e-11$    & $1.68e-10$           & $98$    & 66.10\\
    3            & $4.85e-11$    & $6.36e-10$           & $94$    & 65.09\\
    4            & $1.49e-10$    & $8.55e-10$           & $69$    & 64.62\\ 
    5            & $1.23e-09$    & $7.57e-09$           & $76$    & 78.34\\ 
    6            & $9.09e-10$    & $9.77e-09$           & $67$    & 99.94\\
    7            & $1.71e-09$    & $1.52e-08$           & $77$    & 166.00\\
    8            & $2.80e-09$    & $4.25e-08$           & $80$    & 112.24\\
    \hline
    \end{tabular} 
    \label{tab:16}
\end{table}

Table~\ref{tab:17} presents the impact of different sampling methods on errors, epochs, and computational efficiency. We can see that all three methods achieve good results for this problem. 

\begin{table}[htbp]
    \setlength{\abovecaptionskip}{0cm}
    \setlength{\belowcaptionskip}{0.2cm}
    \centering
    \caption{Impact of sampling methods on errors, epochs, and computational efficiency}
    \renewcommand{\arraystretch}{1.3} 
    \begin{tabular}{ccccc}
    \hline
    Sampling method  & $\|e\|_{L^2}$ & $\|e\|_{L^{\infty}}$ & Epochs  & CPU time(s)\\  \hline
    Uniform          & $2.33e-11$    & $1.68e-10$           & $98$    & 66.10\\
    Gaussian         & $3.42e-11$    & $2.80e-10$           & $119$   & 60.81\\
    Random           & $1.69e-10$    & $1.12e-09$           & $164$   & 74.82\\
    \hline
    \end{tabular}
    \label{tab:17}
\end{table}

We further tested the impact of different random seeds on the results to evaluate the stability and robustness of the DD-SNN method. 
Table~\ref{tab:18} displays the $L^2$ error, $L^\infty$ error, epochs, and computation time under various random seeds. 
The results show that variations in random seeds have minimal impact on accuracy. 
For example, with random seeds set to 0, 202, 2024, and 666666, the $L^2$ errors remain within the $10^{-11}$ range, and the fluctuations in $L^\infty$ errors are also small. 
This indicates that the DD-SNN method can achieve stable accuracy under different random seed conditions, demonstrating its robustness.

\begin{table}[htbp]
    \setlength{\abovecaptionskip}{0cm}
    \setlength{\belowcaptionskip}{0.2cm}
    \centering
    \caption{Impact of random seeds on errors, epochs, and computational efficiency}
    \begin{tabular}{ccccc}
    \hline
    Random seed  & $\|e\|_{L^2}$ & $\|e\|_{L^{\infty}}$ & Epochs        & CPU time(s)\\  \hline
    0            & $4.69e-11$    & $4.06e-10$           & $109$  &230.23\\
    202          & $3.20e-11$    & $2.17e-10$           & $92$   &286.58\\
    2024         & $2.42e-11$    & $5.40e-10$           & $70$   &178.16\\
    666666       & $2.80e-11$    & $2.96e-10$           & $92$   &204.33\\
    888888       & $2.18e-11$    & $2.39e-10$           & $87$   &184.03\\
    \hline
    \end{tabular}
    \label{tab:18}
\end{table}

\begin{figure}
    \centering
    \includegraphics[width=0.8\textwidth]{ 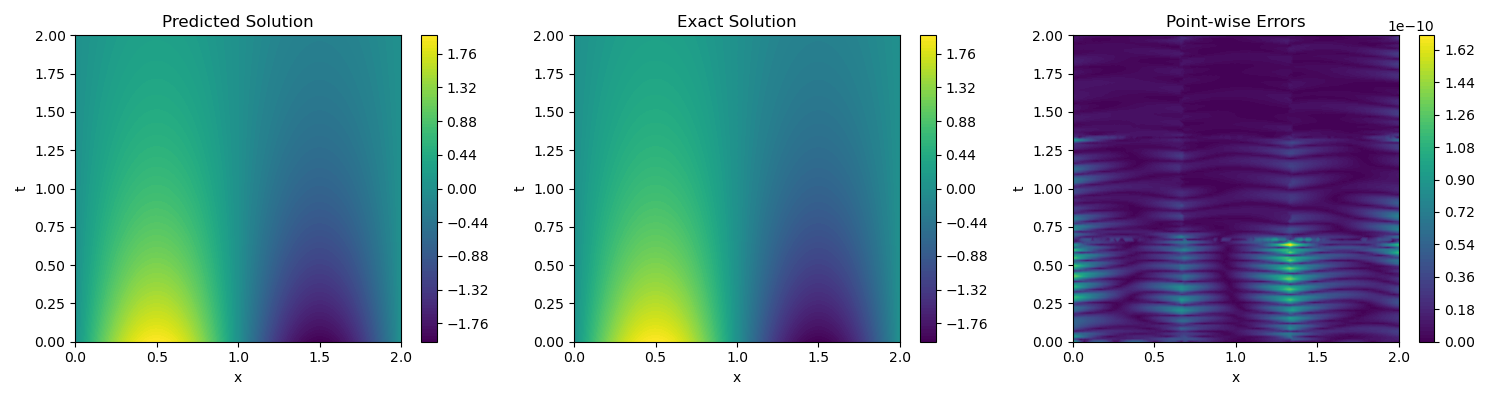}
    \caption{Comparison of numerical and exact solutions for the one-dimensional parabolic equation}
    \label{fig:7}
\end{figure}

Figure~\ref{fig:7} illustrates the global numerical solution and the exact solution of the one-dimensional parabolic equation solved by DD-SNN (refer to Table~\ref{tab:13}), along with the pointwise error distribution. 
The results show that the DD-SNN method can accurately approximate the target solution, verifying its high accuracy and stability for complex problems.

\subsubsection{The Singular Perturbed Boundary Layer Problem}

We consider the singular perturbed boundary layer problem on the domain $\Omega = [a,b]$:
\begin{subequations}
    \begin{align}
        - \varepsilon u_{xx} + u_x &= \varepsilon \pi^2 \sin(\pi x) + \pi \cos(\pi x),\label{eq:39a}\\ 
        u(a) &= h_1,\label{eq:39b}\\
        u(b) &= h_2,\label{eq:39c}
    \end{align}
\end{subequations}
where $u(x)$ is the solution to be determined, and $h_1$ and $h_2$ are constant boundary values. In this experiment, we set $\varepsilon = 0.01$, $a = 0$, $b = 1$, and the exact solution is $u(x)=\sin(\pi x)+\frac{e^{\frac{x}{\varepsilon}}-1}{e^{\frac{1}{\varepsilon}}-1}$.

To ensure fairness, all methods are tested on 800 uniformly sampled points. 
In this experiment, we use the DEEPXDE package to implement the PINNs algorithm, with a fully connected neural network containing 3 hidden layers, each with 100 neurons, and trained for 50,000 epochs. 
The DGM and DRM methods use the same network hyperparameters for consistency. 
The LocELM method divides the solution domain into 8 subdomains, using 3 hidden layers with 100 neurons per layer. 
In the DD-SNN method, the domain is also divided into 8 subdomains (all numerical experiments in this example use this configuration unless otherwise stated), and each subdomain uses a neural network with 2 hidden layers, each with 100 neurons, and a subspace dimension of 100.

Table~\ref{tab:19} presents the errors and epochs for different neural network methods solving Eqs.~\eqref{eq:39a}-\eqref{eq:39c}. 
The results show that after 50,000 epochs, the errors of the PINNs, DGM, and DRM methods remain high, and increasing the number of epochs does not significantly improve accuracy. 
The ELM method produces large errors, while the LocELM method achieves an $L^2$ error of $1.48 \times 10^{-6}$. 
The DD-SNN method achieves an $L^2$ error of $2.45 \times 10^{-7}$ with 63 epochs, demonstrating its advantage in both accuracy and efficiency.

\begin{table}[htbp]
    \setlength{\abovecaptionskip}{0cm}
    \setlength{\belowcaptionskip}{0.2cm}
    \centering
    \caption{Comparison of errors, epochs, and computational efficiency for different methods solving the singular perturbed boundary layer problem}
    \begin{tabular}{ccccc}
    \hline
    Method & $\|e\|_{L^2}$ & $\|e\|_{L^{\infty}}$ & Epochs  & CPU time(s) \\ \hline
    PINN   & $1.26e-03$    & $8.20e-03$           & $50000$ &602.66\\
    DGM    & $1.55e-03$    & $2.95e-03$           & $50000$ &573.82\\
    DRM    & $3.96e-03$    & $1.93e-02$           & $50000$ &572.54\\
    ELM    & $1.40e-01$    & $1.56e-01$           & $-    $ &0.01\\
    LocELM & $1.48e-06$    & $3.94e-06$           & $-    $ &0.14\\
    DD-SNN & $2.45e-07$    & $4.36e-06$           & $63$   &4.59\\ \hline
    \end{tabular}
    \label{tab:19}
\end{table}

\begin{figure}
    \centering
    \includegraphics[width=0.6\textwidth]{ 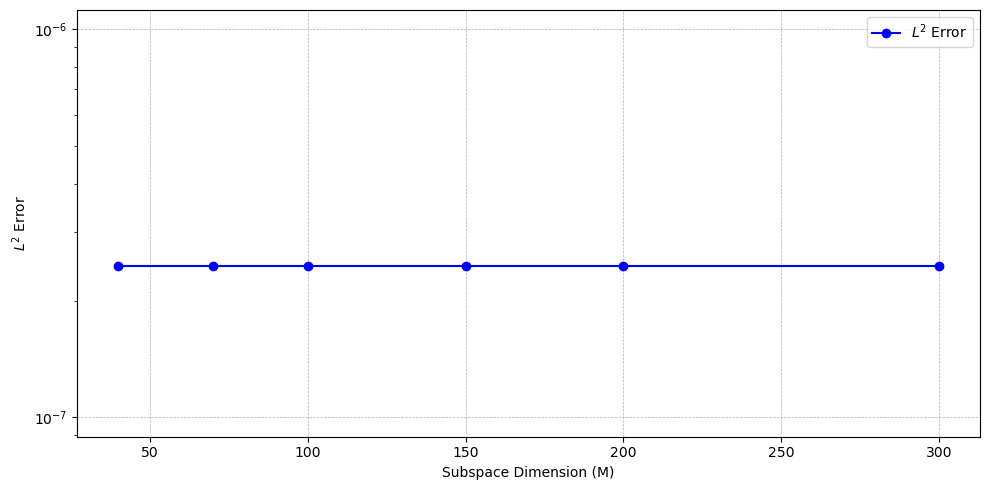}
    \caption{Impact of subspace dimension on $L^2$ error}
    \label{fig:8}
\end{figure}

Figure~\ref{fig:8} shows the subspace dimension plays a limited role in influencing the error at $800$ uniform points, but choosing the appropriate subspace dimension helps to improve the computational efficiency.

Next, we examined the effect of the number of subdomains on errors. 
Keeping the total number of uniform sampling points at 800, the solution domain was divided into varying numbers of subdomains. 
When the number of subdomains is 1, the DD-SNN method degenerates into the standard SNN method. 
Numerical results indicate that for boundary layer problems, increasing the number of subdomains significantly improves solution accuracy. 
When the number of subdomains reaches 6, the error becomes sufficiently small, and further increasing the number of subdomains does not result in noticeable improvements in accuracy.

\begin{table}[htbp]
    \setlength{\abovecaptionskip}{0cm}
    \setlength{\belowcaptionskip}{0.2cm}
    \centering
    \caption{Impact of subdomain numbers on errors and computational efficiency}
    \begin{tabular}{ccccc}
    \hline
    Subdomain    & $\|e\|_{L^2}$ & $\|e\|_{L^{\infty}}$ & Epochs   &CPU time(s)\\ \hline
    1            & $3.49e-01$    & $6.21e-01$           & $72$     & 2.41\\
    2            & $1.40e-02$    & $1.91e-02$           & $62$     & 2.82\\
    4            & $1.35e-06$    & $5.13e-06$           & $49$     & 2.39\\
    6            & $2.47e-07$    & $4.37e-06$           & $52$     & 3.52\\
    8            & $2.45e-07$    & $4.36e-06$           & $63$     & 4.59\\
    16           & $2.23e-07$    & $4.27e-06$           & $65$     & 8.78\\ 
    32           & $2.36e-07$    & $4.38e-06$           & $75$     & 14.00\\\hline
    \end{tabular}
    \label{tab:21}
\end{table}

We further investigated the impact of the number of hidden layers on solution accuracy and computational efficiency, as shown in Table~\ref{tab:22}. 
Keeping the number of subdomains fixed at 8 and each hidden layer containing 100 neurons, the number of hidden layers was varied. The results indicate that the DD-SNN method achieves high accuracy across all tested configurations, with minimal differences in accuracy across different numbers of hidden layers.

\begin{table}[htbp]
    \setlength{\abovecaptionskip}{0cm}
    \setlength{\belowcaptionskip}{0.2cm}
    \centering
    \caption{Impact of hidden layer numbers on errors, epochs, and computational efficiency}
    \begin{tabular}{ccccc}
    \hline
    Hidden layer & $\|e\|_{L^2}$ & $\|e\|_{L^{\infty}}$ & Epochs  &CPU time(s)\\ \hline
    0            & $2.45e-07$    & $4.36e-06$           & $2643$  & 5.69\\
    1            & $2.45e-07$    & $4.36e-06$           & $628$   & 4.25\\
    2            & $2.45e-07$    & $4.36e-06$           & $502$   & 4.59\\
    3            & $2.45e-07$    & $4.36e-06$           & $872$   & 9.29\\
    4            & $2.46e-07$    & $4.38e-06$           & $472$   & 6.74\\ 
    5            & $2.45e-07$    & $4.37e-06$           & $462$   & 7.86\\ 
    6            & $2.45e-07$    & $4.37e-06$           & $1009$  & 16.83\\
    7            & $2.46e-07$    & $4.35e-06$           & $1096$  & 21.00\\
    8            & $2.60e-07$    & $4.31e-06$           & $491$   & 11.56\\
    \hline
    \end{tabular}
    \label{tab:22}
\end{table}

Finally, we tested the effect of three different sampling methods on solution accuracy. 
Numerical experiments were conducted under fixed network hyperparameters (100 neurons per hidden layer, 8 subdomains, 2 hidden layers, and a subspace dimension of 100). 
The results in Table~\ref{tab:23} show that all three sampling methods yield good results.

\begin{table}[htbp]
    \setlength{\abovecaptionskip}{0cm}
    \setlength{\belowcaptionskip}{0.2cm}
    \centering
    \caption{Impact of sampling methods on errors, epochs, and computational efficiency}
    \renewcommand{\arraystretch}{1.3} 
    \begin{tabular}{ccccc}
    \hline
    Sampling method  & $\|e\|_{L^2}$ & $\|e\|_{L^{\infty}}$ & Epochs  & CPU time(s)\\  \hline
    Uniform          & $2.45e-07$    & $4.36e-06$           & $502$   & 4.59\\
    Gaussian         & $5.94e-07$    & $8.58e-06$           & $2141$  & 9.60\\
    Random           & $1.58e-07$    & $1.95e-06$           & $2417$  & 5.45\\
    \hline
    \end{tabular}
    \label{tab:23}
\end{table}

To avoid experimental randomness, we tested the impact of different random seeds on the results. 
Table~\ref{tab:24} displays the effect of various random seeds on errors, epochs, and computational efficiency.

\begin{table}[htbp]
    \setlength{\abovecaptionskip}{0cm}
    \setlength{\belowcaptionskip}{0.2cm}
    \centering
    \caption{Impact of random seeds on errors, epochs, and computational efficiency}
    \begin{tabular}{ccccc}
    \hline
    Random seed  & $\|e\|_{L^2}$ & $\|e\|_{L^{\infty}}$ & Epochs  & CPU time(s)\\  \hline
    1            & $2.45e-07$    & $4.36e-06$           & $3399$  & 7.39\\
    202          & $2.45e-07$    & $4.36e-06$           & $502$   & 4.59\\
    666          & $2.45e-07$    & $4.36e-06$           & $3407$  & 7.09\\
    888          & $2.45e-07$    & $4.36e-06$           & $2578$  & 5.66\\
    2024         & $2.45e-07$    & $4.36e-06$           & $2619$  & 5.93\\
    666888       & $2.45e-07$    & $4.36e-06$           & $3325$  & 7.29\\
    \hline
    \end{tabular}
    \label{tab:24}
\end{table}

\begin{figure}
    \centering
    \includegraphics[width=0.8\textwidth]{ 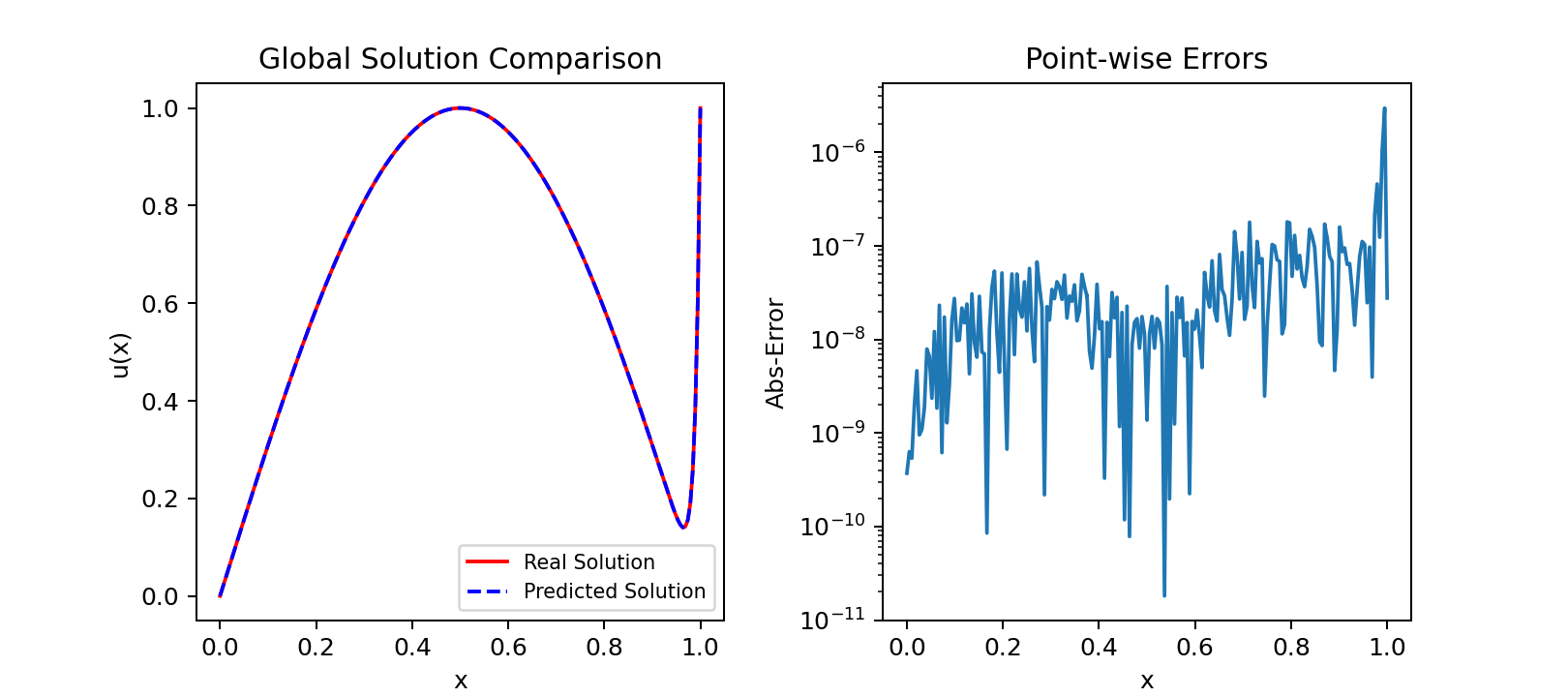}
    \caption{Comparison of numerical and exact solutions for the one-dimensional convection-diffusion equation}
    \label{fig:9}
\end{figure}

Figure~\ref{fig:9} illustrates the comparison of numerical and exact solutions for the singular perturbed boundary layer problem solved by DD-SNN, along with the pointwise error distribution. 
The results show that DD-SNN accurately approximates the target solution, validating its high accuracy and stability for this class of problems.

\subsection{Nonlinear Example}

\subsubsection{One-Dimensional Nonlinear Helmholtz Equation}
We evaluate the performance of the DD-SNN method using the boundary value problem of a one-dimensional nonlinear Helmholtz equation. 
Consider the boundary value problem on the domain $\Omega = [a,b]$:
\begin{subequations}
    \begin{align}
        &\frac{\mathrm{d}^2 u}{\mathrm{d}x^2} - \lambda u + c\sin(u) = f(x), \label{eq:40a} \\
        &u(a) = h_1,  \label{eq:40b} \\
        &u(b) = h_2,  \label{eq:40c}
    \end{align}
    \label{eq:40}
\end{subequations}
where $u(x)$ is the solution to be determined, $f(x)$ is the given source term, $\lambda$ and $c$ are constant parameters, and $h_1$ and $h_2$ are boundary values. 
For the parameters in the equations and the domain, we set:
\begin{equation}
    \lambda = 50, \quad c = 10, \quad a = 0, \quad b = 8.
\end{equation}
The source term $f(x)$ and boundary values $h_1$ and $h_2$ are chosen so that the following function satisfies Eq.~\eqref{eq:40}:
\begin{equation}
    u(x) = \sin\left(3\pi x+\frac{3\pi}{20}\right)\cos\left(4\pi x-\frac{2\pi}{5}\right) +\frac{3}{2}+\frac{x}{10}.
\end{equation}

Table~\ref{tab:25} presents the errors and epochs for different neural network methods when solving this equation. 
To ensure fairness, all methods are tested on 800 uniformly sampled points. We use the PINNs algorithm in the DEEPXDE package, with a network structure of 3 hidden layers, where the first two layers have 100 neurons each and the last layer has 200 neurons, trained for 50,000 epochs. 
DGM and DRM use the same network hyperparameters. 
For LocELM, nonlinear least squares is used with a convergence tolerance of $10^{-6}$. 
The results show that even after 50,000 epochs, the errors of PINNs, DGM, and DRM remain high, and increasing the number of epochs does not significantly improve accuracy.

The DD-SNN method divides the solution domain into 8 subdomains (all experiments in this example use 8 subdomains), employs a network with two hidden layers, each containing 100 neurons, a subspace dimension of 200, and sets the maximum number of Picard iterations to 20 with a convergence tolerance of $10^{-6}$. 
Under these settings, DD-SNN achieves significantly fewer total iterations and substantially higher accuracy than other methods, demonstrating its advantages in both accuracy and computational cost.

\begin{table}[htbp]
    \setlength{\abovecaptionskip}{0cm}
    \setlength{\belowcaptionskip}{0.2cm}
    \centering
    \caption{Comparison of errors, epochs, and computational efficiency for solving the nonlinear Helmholtz equation}
    \begin{tabular}{cccccc}
    \hline
    Method & $\|e\|_{L^2}$ & $\|e\|_{L^{\infty}}$ & Epochs    &Picard Iterations & CPU time(s) \\ \hline
    PINN   & $7.11e-01$    & $2.41e+00$           & $50000$   & $-$      &1596.26\\
    DGM    & $6.27e-03$    & $1.65e-02$           & $50000$   & $-$      &343.34\\
    DRM    & $1.29e-02$    & $3.92e-02$           & $50000$   & $-$      &619.43\\
    LocELM & $8.51e-07$    & $3.92e-06$           & $-    $   & $-$      &179.88\\
    DD-SNN & $9.49e-10$    & $3.87e-09$           & $2357$    & $8$      &131.22\\ \hline
    \end{tabular}
    \label{tab:25}
\end{table}

\begin{figure}
    \centering
    \includegraphics[width=0.6\textwidth]{ 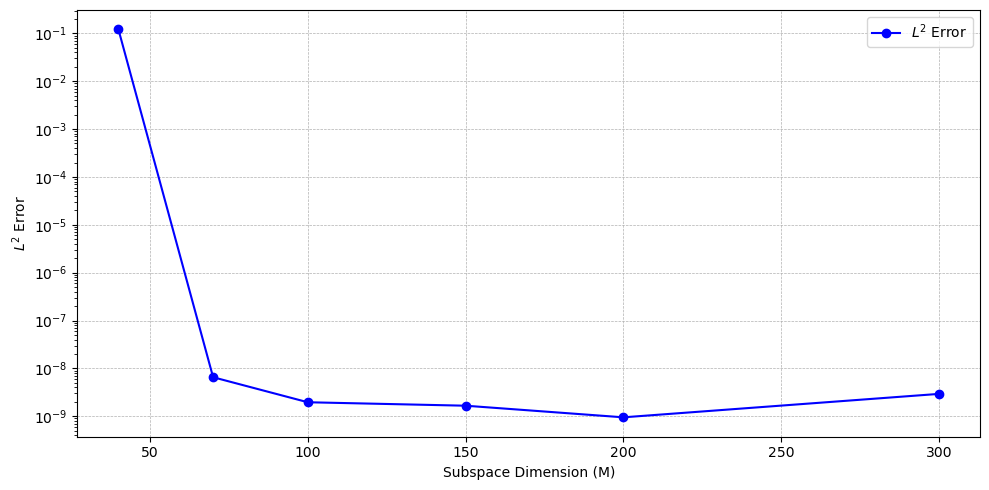}
    \caption{Impact of subspace dimension on $L^2$ error}
    \label{fig:10}
\end{figure}

Figure~\ref{fig:10} shows the impact of different subspace dimensions ($M$) on error, epochs, and nonlinear iterations at $800$ uniform points. 
The convergence tolerance for Picard iterations is set to $10^{-6}$, and convergence is typically achieved in approximately 8 iterations. 
The results indicate that as the subspace dimension $M$ increases, the $L^2$ error decreases. 
This trend highlights the importance of higher subspace dimensions for improving solution accuracy. 
Although increasing subspace dimensions slightly increases epochs, the impact on nonlinear iterations is minimal, enabling DD-SNN to maintain a balance between high accuracy and computational cost.

Table~\ref{tab:27} details the impact of different numbers of hidden layers on errors, epochs, and computational time. 
The experiments fixed the subspace dimension to 200, used 800 uniformly sampled points, and divided the domain into 8 subdomains to explore how changes in the number of hidden layers affect the results. 
When the number of hidden layers is 2, the $L^2$ error reaches its minimum value of $9.49 \times 10^{-10}$, and the $L^\infty$ error is $3.87 \times 10^{-9}$, achieving the best accuracy. 
A two-layer hidden structure offers the best balance between accuracy and computational cost, providing high precision with a moderate number of nonlinear iterations.

\begin{table}[htbp]
    \setlength{\abovecaptionskip}{0cm}
    \setlength{\belowcaptionskip}{0.2cm}
    \centering
    \caption{Impact of hidden layer numbers on errors, epochs, and computational efficiency}
    \begin{tabular}{cccccc}
    \hline
    Hidden layer& $\|e\|_{L^2}$ & $\|e\|_{L^{\infty}}$ & Epochs   &Picard Iterations & CPU time(s)\\ \hline
    0            & $6.49e-08$    & $4.73e-07$           & $5000$  &8    & 99.62\\
    1            & $1.06e-08$    & $4.96e-08$           & $2140$  &8    & 91.40\\
    2            & $9.49e-10$    & $3.87e-09$           & $2357$  &8    & 131.22\\
    3            & $1.25e-09$    & $5.50e-09$           & $1665$  &8    & 160.32\\
    4            & $2.17e-09$    & $8.67e-09$           & $1943$  &8    & 184.79\\ 
    5            & $3.85e-09$    & $1.91e-08$           & $1526$  &8    & 239.08\\ 
    6            & $5.59e-09$    & $2.41e-08$           & $1019$  &8    & 199.04\\
    7            & $7.57e-09$    & $3.53e-08$           & $631$   &8    & 92.30\\
    8            & $1.58e-08$    & $9.11e-08$           & $393$   &8    & 98.37\\
    \hline
    \end{tabular}
    \label{tab:27}
\end{table}

Table~\ref{tab:28} shows the effect of the number of subdomains on errors, iterations, and CPU time. 
As the number of subdomains increases to 32, the $L^2$ error significantly decreases from $1.04 \times 10^0$ to $2.56 \times 10^{-11}$, and the number of epochs also reduces. 
For instance, 32 subdomains require only 77 epochs and 94.27 seconds of computation time.

\begin{table}[htbp]
    \setlength{\abovecaptionskip}{0cm}
    \setlength{\belowcaptionskip}{0.2cm}
    \centering
    \caption{Impact of subdomain numbers on errors and computational efficiency}
    \begin{tabular}{cccccc}
    \hline
    Subdomain    & $\|e\|_{L^2}$ & $\|e\|_{L^{\infty}}$ & Epochs   &Picard Iterations &CPU time(s)\\ \hline
    1            & $1.49e+01$    & $1.38e+02$           & $5000$  &20   & 86.72\\
    2            & $1.04e+00$    & $7.39e+00$           & $5000$  &20   & 169.70\\
    4            & $4.56e-08$    & $2.93e-07$           & $5000$  &8    &266.89\\
    6            & $1.67e-09$    & $1.15e-08$           & $2478$  &8    &236.87\\
    8            & $9.49e-10$    & $3.87e-09$           & $2357$  &8    &131.22\\
    16           & $1.74e-10$    & $8.98e-10$           & $396$   &8    & 117.19\\ 
    32           & $2.56e-11$    & $6.97e-11$           & $77$    &8    &94.27\\ \hline
    \end{tabular}
    \label{tab:28}
\end{table}

Table~\ref{tab:29} presents the impact of different sampling methods on errors, epochs, and computational efficiency. 
The experiments demonstrate that using Gaussian sampling achieves the best performance, with an $L^2$ error of $1.53 \times 10^{-10}$ and $1907$ epochs. 
Uniform sampling is almost as good as Gaussian sampling.
In contrast, random sampling is slightly less accurate compared to the first two, but it also achieves a high level of accuracy with almost constant computational efficiency.

\begin{table}[htbp]
    \setlength{\abovecaptionskip}{0cm}
    \setlength{\belowcaptionskip}{0.2cm}
    \centering
    \caption{Impact of sampling methods on errors, training epochs, and computational efficiency}
    \renewcommand{\arraystretch}{1.3} 
    \begin{tabular}{cccccc}
    \hline
    Sampling method & $\|e\|_{L^2}$ & $\|e\|_{L^{\infty}}$ & Epochs  &Picard Iterations & CPU time(s)\\  \hline
    Uniform         & $9.49e-10$    & $3.87e-09$           & $2357$  & $8$ &122.43\\
    Gaussian        & $1.53e-10$    & $7.82e-10$           & $1907$  & $8$ &134.75\\
    Random          & $4.43e-08$    & $1.60e-07$           & $2357$  & $8$ &125.50\\
    \hline
    \end{tabular}
    \label{tab:29}
\end{table}

\begin{figure}
    \centering
    \includegraphics[width=0.8\textwidth]{ 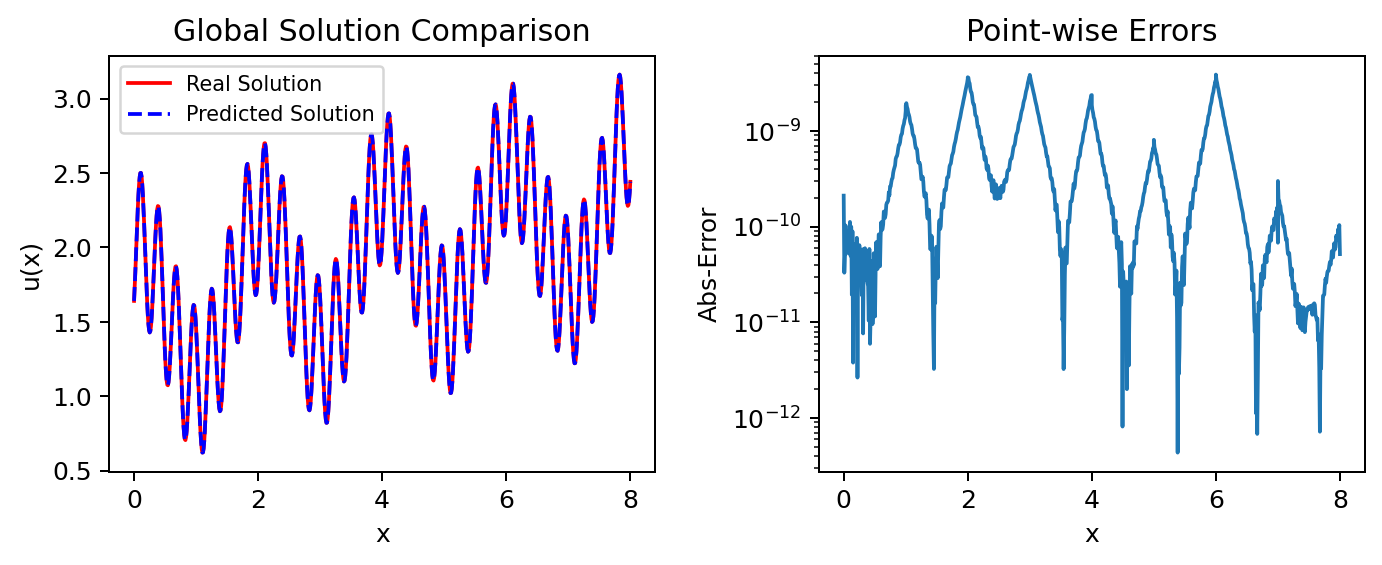}
    \caption{Error distribution for the one-dimensional nonlinear Helmholtz equation}
    \label{fig:11}
\end{figure}

Figure~\ref{fig:11} illustrates the global numerical solution and the exact solution for the one-dimensional nonlinear Helmholtz equation solved by DD-SNN (refer to Table~\ref{tab:25}), along with the pointwise error distribution. 
The results show that DD-SNN accurately approximates the target solution while maintaining minimal pointwise errors, validating its high accuracy and stability for nonlinear problems.

\subsubsection{One-Dimensional Viscous Burgers Equation}

We consider the one-dimensional viscous Burgers equation on the spatiotemporal domain $\Omega = [a,b] \times [0,T]$:
\begin{subequations}
    \begin{align}
        \frac{\partial u}{\partial t} + u\frac{\partial u}{\partial x} &= \nu \frac{\partial^2 u}{\partial x^2} + f(x,t), \label{eq:43a} \\ 
        u(a,t) &= g_1(t), \label{eq:43b} \\ 
        u(b,t) &= g_2(t), \label{eq:43c} \\ 
        u(x,0) &= h(x), \label{eq:43d}
    \end{align}
\end{subequations}
where $u(x,t)$ is the solution to be determined, $\nu$ is the viscosity coefficient, $f(x,t)$ is the source term, and $g_1(t)$ and $g_2(t)$ are the time-dependent boundary conditions at $x=a$ and $x=b$, respectively. 
The initial condition $h(x)$ describes the state of $u(x,t)$ at $t=0$. For this problem, we set:
\begin{equation}
    \nu = 0.01, \quad a = 0, \quad b = 2 \pi, \quad T = 1.
\end{equation}
The source term $f(x,t)$ and boundary values $g_1(t)$ and $g_2(t)$ are chosen such that the following function satisfies Eqs.~\eqref{eq:43a}-\eqref{eq:43d}:
\begin{equation}
    u(x,t) = e^{-0.01t}\sin(x).
    \label{eq:45}
\end{equation}

To ensure fairness, all methods were tested on $20 \times 20 \times 8$ uniform sampling points. 
Table~\ref{tab:30} compares the performance of PINNs, DGM, DRM, and DD-SNN methods for solving the one-dimensional viscous Burgers equation. 
After 50,000 epochs, the $L^2$ errors for PINNs, DGM, and DRM methods were $2.63 \times 10^{-1}$, $7.39 \times 10^{-1}$, and $2.30 \times 10^{-1}$, respectively, with computational times of 1425.60 seconds, 4003.62 seconds, and 634.24 seconds. 
These results highlight the limitations of these methods in reducing errors and achieving computational efficiency.

In contrast, the DD-SNN method divides the solution domain into 8 subdomains, using a network with two hidden layers, 100 neurons per layer, and a subspace dimension of 200. 
The Picard iteration convergence tolerance $\varepsilon_{non}$ is set to $1 \times 10^{-12}$, and Newton iteration $\varepsilon_{non}$ is set to $1 \times 10^{-13}$. 
Each subdomain requires an average of only 94 epochs. 
Using Picard iterations, the $L^2$ error was reduced to $3.56 \times 10^{-10}$ and the $L^{\infty}$ error to $2.10 \times 10^{-9}$ within 17 iterations. 
Switching to Newton iterations, the $L^2$ error was further reduced to $4.19 \times 10^{-10}$ after just 2 iterations. 
This demonstrates the superior accuracy and efficiency of the DD-SNN method. 
LocELM, using the NLSQ-perturn strategy \cite{Dong2021LocELM}, failed to converge after 2 days, highlighting the computational advantages of DD-SNN.

\begin{table}[htbp]
    \setlength{\abovecaptionskip}{0cm}
    \setlength{\belowcaptionskip}{0.2cm}
    \centering
    \caption{Comparison of errors, epochs, and computational efficiency for solving the one-dimensional viscous Burgers equation}
    \begin{tabular}{cccccc}
    \hline
    Method           & $\|e\|_{L^2}$ & $\|e\|_{L^{\infty}}$ & Epochs & Iterations & CPU time(s)\\ \hline
    PINN             & $2.63e-01$    & $5.74e-01$           & $50000$ & $-$       & 1425.60\\
    DGM              & $7.39e-01$    & $1.41e+00$           & $50000$ & $-$       & 4003.62\\
    DRM              & $2.30e-01$    & $3.16e-01$           & $50000$ & $-$       & 634.24\\
    LocELM           & $-$           & $-$                  & $-$     & $-$       & $-$\\
    DD-SNN (Picard)  & $3.56e-10$    & $2.10e-09$           & $94$    & $17$      & 23.56\\ 
    DD-SNN (Newton)  & $4.19e-10$    & $2.46e-09$           & $94$    & $2$       & 24.89\\ \hline
    \end{tabular}
    \label{tab:30}
\end{table}

\begin{figure}
    \centering
    \includegraphics[width=0.6\textwidth]{ 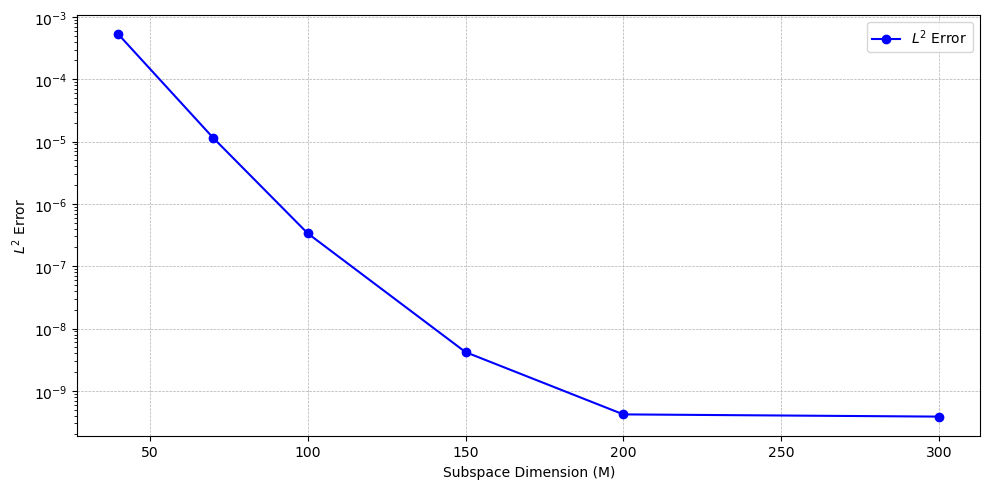}
    \caption{Impact of subspace dimension on $L^2$ error}
    \label{fig:12}
\end{figure}

Figure~\ref{fig:12} shows the changes in $L^2$ error and epochs under different subspace dimensions at $20\times20\times8$ uniform points. 
The experiments used a network with two hidden layers, 100 neurons per layer, and divided the solution domain into 8 subdomains. 
The results indicate that as the subspace dimension $M$ increases, the $L^2$ error decreases significantly.

Table~\ref{tab:32} analyzes the impact of the number of subdomains on $L^2$ error, $L^{\infty}$ error, epochs, and CPU time, using a subspace dimension of $200$ and $20 \times 20$ sampling points per subdomain. 
Numerical results indicate that increasing the subdomain count from $1 \times 1$ to $2 \times 2$ reduces the $L^2$ error from $1.39 \times 10^{-7}$ to $1.83 \times 10^{-9}$, with a similar decrease in $L^{\infty}$ error. 
The configuration of $4 \times 2$ subdomains achieves high precision while controlling computational costs, making it the optimal balance.

\begin{table}[htbp]
    \setlength{\abovecaptionskip}{0cm}
    \setlength{\belowcaptionskip}{0.2cm}
    \centering
    \caption{Impact of subdomain numbers on errors and training efficiency}
    \begin{tabular}{cccccc}
    \hline
    Subdomain   & $\|e\|_{L^2}$ & $\|e\|_{L^{\infty}}$  & Epochs    & Iterations & CPU time(s)\\ \hline
    $1\times1$   & $1.39e-07$    & $6.17e-07$           & $5000$     &$20$      & 90.23\\
    $2\times2$   & $1.83e-09$    & $7.99e-09$           & $1235$     &$20$      & 191.09\\
    $4\times2$   & $4.19e-10$    & $2.46e-09$           & $94$       &$2$       & 24.89\\
    $4\times4$   & $9.35e-10$    & $5.86e-09$           & $143$      &$3$       & 73.00\\ 
    $6\times4$   & $7.10e-10$    & $3.43e-09$           & $132$      &$3$       & 159.64\\ \hline
    \end{tabular}
    \label{tab:32}
\end{table}

Table~\ref{tab:33} examines the influence of the number of hidden layers on $L^2$ error, $L^{\infty}$ error, epochs, and CPU time, using a subspace dimension of $M=200$ and $4 \times 2$ subdomains. 
The results demonstrate that all tested configurations achieve convergence, but increasing the number of hidden layers results in higher training costs. 
When the number of hidden layers is 2, the $L^2$ error is $4.19 \times 10^{-10}$, achieving a good balance between accuracy and computational efficiency.

\begin{table}[htbp]
    \setlength{\abovecaptionskip}{0cm}
    \setlength{\belowcaptionskip}{0.2cm}
    \centering
    \caption{Impact of hidden layer numbers on errors, epochs, and computational efficiency}
    \begin{tabular}{cccccc}
    \hline
    Hidden layer& $\|e\|_{L^2}$  & $\|e\|_{L^{\infty}}$ & Epochs   & Iterations  & CPU time(s)\\ \hline
    0            & $1.13e-10$    & $5.97e-10$           & $353$    &$2$      & 24.86\\
    1            & $1.76e-10$    & $1.07e-09$           & $86 $    &$3$      & 16.48\\
    2            & $4.19e-10$    & $2.46e-09$           & $94$     &$2$      & 24.89\\
    3            & $7.33e-10$    & $4.92e-09$           & $313$    &$10$     & 77.52\\
    4            & $1.36e-09$    & $8.96e-09$           & $662$    &$20$     & 182.65\\ 
    5            & $2.71e-09$    & $1.68e-08$           & $328$    &$20$     & 124.49\\ 
    6            & $5.32e-09$    & $2.92e-08$           & $566$    &$20$     & 209.40\\
    7            & $2.27e-08$    & $1.07e-07$           & $1137$   &$20$     & 423.81\\
    8            & $3.58e-08$    & $2.18e-07$           & $528$    &$20$     & 241.10\\
    \hline
    \end{tabular}
    \label{tab:33}
\end{table}

\begin{figure}
    \centering
    \includegraphics[width=0.8\textwidth]{ 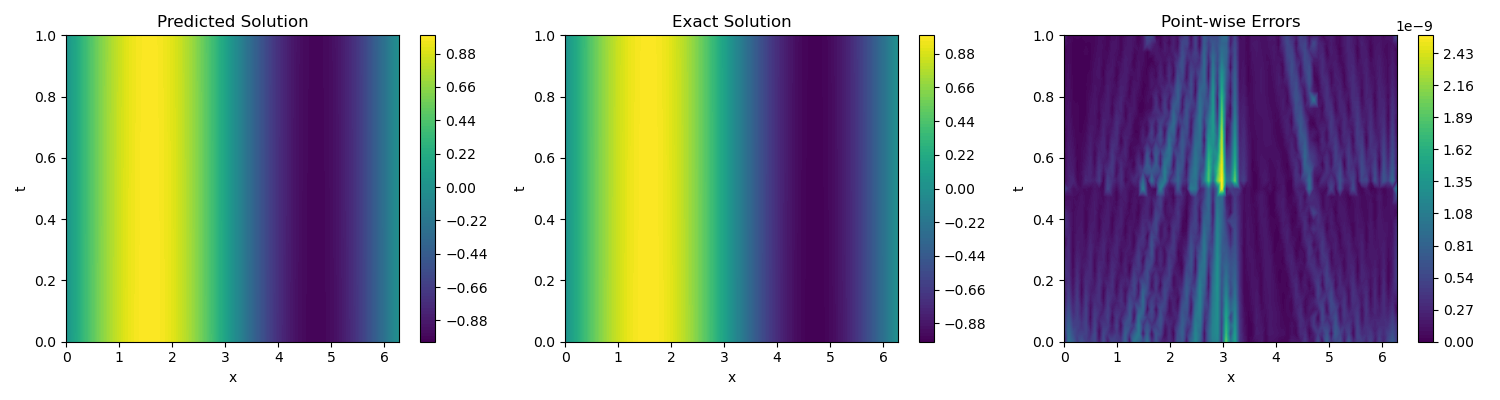}
    \caption{Comparison of numerical and exact solutions for the one-dimensional viscous Burgers equation}
    \label{fig:13}
\end{figure}

Figure~\ref{fig:13} illustrates the numerical solution obtained by the DD-SNN (Newton) method for the one-dimensional viscous Burgers equation, as presented in Table~\ref{tab:30}. 
The comparison between the numerical solution and the exact solution shows that the DD-SNN method can accurately approximate the target solution, with minimal pointwise error distribution, demonstrating its high accuracy and stability for nonlinear problems.

\section{Conclusion}
This paper proposed the domain decomposition subspace neural network method for efficiently solving linear and nonlinear partial differential equations. 
By combining domain decomposition and subspace neural networks, our method achieved high accuracy with low computational cost. 
Our method solved various problems, including the linear and nonliner Helmholtz equation, diffusion equation, singular perturbed boundary layer problem, and viscous Burgers equation. 
Numerical experiments demonstrated that DD-SNN significantly outperformed existing methods, such as PINNs, DGM, DRM and ELM, highlighting its potential for solving complex PDEs.

\section*{Acknowledgments}
This work is partially supported by National Natural Science Foundation of China (12071045) and Fund of National Key Laboratory of Computational Physics.

\end{document}